\documentclass[a4paper,10pt]{article}
\usepackage{amsthm}
\usepackage[utf8]{inputenc}   %paquet que serveix per poder escriure
\usepackage[english]{babel}
\usepackage{amsfonts}
\usepackage{amssymb}
\usepackage{amsmath}
\usepackage{doc}
\usepackage{fontenc}
\usepackage{ifthen}
\usepackage{latexsym}
\usepackage{makeidx}
\usepackage{graphicx}
\usepackage{float}
\usepackage{array}
\usepackage[all]{xy}
\usepackage{color}
\usepackage[hyperref]{degt}

\def\2{\color{red}}
\def\3{\color{magenta}}
\def\4{\color{cyan}}
\def\5{\color{blue}}
\let\sectionname\sectionautorefname
\let\subsectionname\subsectionautorefname

\theoremstyle{remark}

\newcommand{\beq}{\begin{equation} }
\newcommand{\enq}{\end{equation}}

\let\ra\undefined
\newcommand{\ra}{\rightarrow}

\DeclareMathOperator{\sing}{\rm sing}
\newcommand{\KK}{{\mathbb K}}

\newcommand{\XXf}{{S}}
\def\tX{\tilde{\XXf}}
\def\KK{\tilde{K}}
\def\hh{h}
\newcommand{\LLL}{\ell}

\newcommand{\sLLL}{\ell}
\def\cross{_{{\!\times\!}}}
\def\LLX{\ell\cross}
\def\SSf{\Sigma}

\def\GQ{\operatorname{GQ}}

\def\I{\mathrm{I}}
\def\II{\mathrm{II}}
\def\III{\mathrm{III}}
\def\IV{\mathrm{IV}}
\def\chtop{\chi_\mathrm{top}}

\def\bX{\bold X}
\def\bY{\bold Y}
\def\bZ{\bold Z}
\def\bW{\bold W}
\def\bT{\bold T}
\def\bS{\bold S}
\def\bQ{\bold Q}
\def\bP{\bold P}
\def\bH{\bold H}

\def\FF{\mathbb{F}}
\def\bu{\bold{u}}
\def\bv{\bold{v}}
\def\one{\bold{c}}
\def\pp{q}

\def\mutab#1{\vtop\bgroup\offinterlineskip\def\*{\rlap{$^*$}}%
	\ialign\bgroup\vrule height11pt depth3.5pt\relax
	\quad\hss$##$\hss\quad\vrule&\quad\hss$##$\hss&
	\quad\hss$##$\hss\quad&\vrule\quad\hss$##$\hss\quad\vrule\cr
	\noalign{\hrule height0pt}\multispan4\strut\hss#1\hss\cr
	\noalign{\smallskip\hrule}
	S&\omit\span\omit\hss Elliptic fiber\hss&\mu\cr
	\noalign{\hrule}}
\def\endmutab{\crcr
	\noalign{\hrule}\crcr\egroup\egroup}

{
\catcode`\1\active
\catcode`\.\active
\gdef\tconfig{\vcenter\bgroup\offinterlineskip
\catcode`\1\active\catcode`\.\active\def1{&\bullet}\def.{&\cdot}%
\halign\bgroup
 \hbox to0pt{\hss##\hss}&&\hbox to9pt{\small\hss$\mathstrut##$\hss}\cr}
\gdef\endtconfig{\crcr\egroup\egroup}
}

\def\bJ{\bold J}
\def\bE{\bold E}
\def\tA{\tilde{\bold A}}
\def\tD{\tilde{\bold D}}
\def\tE{\tilde{\bold E}}

\def\conic{C}

\def\qqser#1{(\mathrm{Q}#1)}
\pdef\qser#1{\ifmmode\qqser{#1}\else$\qqser{#1}$\fi}
\def\AAser#1{\ifcase#1\relax\or\bold\conic\or\bP\or\bT\or\bX\or\bJ\*\or\bJ\fi}
\def\Aser#1{\ifmmode\AAser{#1}\else$\AAser{#1}$\fi}
\let\ser\Aser

\def\minifirst#1{\hss$#1$}
\def\minitab#1{\vcenter\bgroup\rm
 \def\-##1{\setbox0\hbox{$00$}\hbox to\wd0{\hss$##1$\hss}}%
 \def\.{\cdot}
 \let\\\cr
 \halign\bgroup\strut\minifirst{##}&&#1\hss$##$\hss\cr}
\def\endminitab{\crcr\egroup\egroup}

\def\HHH{\setbox0\hbox{$(0)$}\hbox to\wd0{\hss$\HH$\hss}}

\let\ul\underline

\def\*{^\star}
\def\rat{_{\mathrm{rat}}}

\def\Fn{\operatorname{Fn}}
\def\idZ{_{-2}}
\def\idK{_{-1}}
\def\FnZ{\Fn\idZ}
\def\FnK{\Fn\idK}
\def\qZ{q\idZ}
\def\qK{q\idK}
\def\Vset{\operatorname{V}}
\def\VZ{\Vset\idZ}
\def\VK{\Vset\idK}
\def\Eset{\operatorname{E}}
\def\Eall{\operatorname{\tilde E}}
\def\Emax{\Eall_{\max}}
\def\Cset{\operatorname{C}}
\def\Vec{\operatorname{vec}}
\def\bnd{\operatorname{bnd}}
\def\VecZ{\Vec^+}
\def\VecK{\Vec}
\def\bndZ{\bnd^+}
\def\bndK{\bnd}

\def\dual{^\vee}

\def\be{\bar{e}}

\def\vv#1{\bold1_{#1}}
\def\IS{\Cal I}

\def\ssing{\Delta}

\title{On large configurations of lines on quartic surfaces}

\author{Alex Degtyarev\thanks{Partially supported by the T\"{U}B\.ITAK grant
123F111}\ \
and S\l awomir Rams\thanks{Research funded by
		the National Science Centre, Poland, Opus  grant
		no.\ 2024/\allowbreak 53/\allowbreak B/\allowbreak ST1/\allowbreak 01413}}

\date{}

\numberwithin{equation}{section}

\begin{document}

\let\sectionautorefname\sectionname
\let\subsectionautorefname\subsectionname

\maketitle
\begin{abstract}
We estimate the number of lines on a non-$K3$ quartic surface.
Such a surface with only isolated double point(s)
contains at most twenty lines;
this bound is attained
by a unique configuration
of lines
and
by a surface with a certain limited set of singularities.
%only when the surface
%in question has one non-simple {\5 singularity of the type $\bX_{1,0}$,  $\bX_{1,1}$ or $\bZ_{11}$} \mnote{s:to changed}
% (Theorem~\ref{th.X}).
We have similar
itemized bounds for other types of non-simple singularities, which
culminate in at most $31$ lines
on a non-$K3$ quartic not ruled by lines;
this bound is only attained on the quartic monoids described
by K.~Rohn.
\end{abstract}

%%%%%%%%%%%%%%%%%%%%%%%%%%%%%%%%%%%%%%%%%%%%%%%%%%%%%%%%%%%%%%%%%%%%
%%%%%%%%%%%%%%%%%%%%%%%%%%%%%%%%%%%%%%%%%%%%%%%%%%%%%%%%%%%%%%%%%%%%
%%%%%%%  Sec Introduction
%%%%%%%%%%%%%%%%%%%%%%%%%%%%%%%%%%%%%%%%%%%%%%%%%%%%%%%%%%%%%%%%%%%%

\section{Introduction} \label{S.introduction}

Line configurations on surfaces in  $ \PP^3(\CC)$
have been studied since the
XIXth century.
For example, the
maximal
number of lines
on a cubic with a given
set of singularities
is known (see, \eg, \cite[Table~1]{sakamaki} and the bibliography therein).
In contrast, far less is known for surfaces of higher degree.
Line configurations on non-normal
quartic surfaces, resp.\ quartics with an isolated triple point, were
studied, among others, by Clebsch, Cremona, Segre, Rohn (see, \eg, \cite{Clebsch-1868,top-11}
and the bibliography therein), resp.\ Rohn~\cite{rohn-dreifache-punkt}.
Partial results on other quartics can be found in~\cite{jessop}, but
most questions on line configurations on $K3$-quartics
%have been answered within
were not answered until
the last decade (\cite{RS,dis-2015,veniani-2013,DR-2023}),
following
the seminal paper~\cite{Segre}.
Ultimately, the maximal number is $64$ in the smooth case, see
\cite{Segre,RS,dis-2015}, and it drops down to $52$ in the presence of at
least one simple singular point, see~\cite{DR-2023}.
Some bounds on the number of lines on \emph{non-$K3$}
quartics
(the principal subject of the present paper)
can be found in~\cite{gonzalezalonso}, where a sharp upper bound for
affine complex quartics is
obtained.

After~\cite{gonzalezalonso,DR-2023} it is generally understood that,
roughly, the more complex the singularities of a quartic are the fewer
lines it may have.
Thus, the main aim of the present paper is to reconfirm this observation by
completing the few missing cases and finding
upper
bounds on the cardinality
of the configurations of lines
on complex projective (necessarily irreducible)
non-$K3$
quartic
surfaces with various types of singularities
\emph{provided that the surface is not ruled by lines}.
To avoid the ambiguity in the case of a line of
singular points, we agree that a \emph{line} is
a degree-one curve in $\PP^3(\CC)$.

In particular, we prove the following
statements.
(Throughout the paper, we use the classification of isolated
hypersurface singularities found in \cite{AVG1}.)

\theorem[see \autoref{proof.main}]\label{thm-non-K3-projective-double}
Let
$\XXf \subset \PP^3(\CC)$
be a non-$K3$ quartic surface
with at worst isolated double points as singularities.
Then
%the quartic
$\XXf$ contains at most $20$ lines, and this bound is sharp.
\endtheorem

%\addendum[see \autoref{proof-add-non-K3-projective-double}]\label{add-non-K3-projective-double}
%A quartic $\XXf$ as in \autoref{thm-non-K3-projective-double} has
%either $20$ or at most $18$ lines. In the former case, $\XXf$ has a
%%unique non-simple singularity of the
%{\3single singular point~$O$, which is of}
%type
%$\bX_{1,0}$,  $\bX_{1,1}$ or $\bZ_{11}$, and the $20$ lines
%consist of four coplanar lines {\3in the tangent cone $C_O\XXf$}
%and a generalized quadrangle $\GQ(3,1)$
%\rom(see \autoref{th.X} and \autoref{rem.X} for a more precise
%description\rom).
%%Moreover,\mnote{\3Th too loaded}
%%if $\XXf$ contains $20$ lines, then it has a
%%%unique non-simple singularity of the
%%{\3single singular point~$O$, which is of}
%%type
%%$\bX_{1,0}$,  $\bX_{1,1}$ or $\bZ_{11}$, and the $20$ lines
%%\mnote{\3 probably, we should speculate about what we do or do not know later,
%%in \autoref{rem.X}? Including the realizability of each of the three types}
%%consist of four coplanar lines {\3in the tangent cone of~$\XXf$ at~$O$}
%%and a generalized quadrangle $\GQ(3,1)$
%%\rom(see \autoref{th.X} and \autoref{rem.X} below for a precise description\rom).
%%%of the configuration of $20$ lines see  Theorem~\ref{th.X} below\rom).
%\endaddendum

\addendum[see \autoref{proof-add-non-K3-projective-double}]\label{add-non-K3-projective-double}
A quartic $\XXf\subset\Cp3(\C)$ as in \autoref{thm-non-K3-projective-double} has
either $20$ or at most $18$ lines.
Furthermore\rom:
\roster
\item\label{add.8}
if $\XXf$ has more than $8$ lines, then it is rational\rom;
in particular,
$\XXf$ has a single non-simple singular point~$O$ \rom(\cf. \autoref{lem.rational}\rom)\rom;
\item\label{add.16}
if $\XXf$ has more than $16$ lines, then this point $O$ is adjacent to
$\bX_9:=\bX_{1,0}$
\rom(see the \ser4-series in \autoref{conv.types}\rom)\rom;
\item\label{add.18}
if $\XXf$ has $18$ lines, then
%$O$ is of type
%$\bX_{1,0}$,  $\bX_{1,1}$, or $\bZ_{11}$\rom; in addition, the set $\sing(\XXf)$ may
%have at most one simple node $\bA_1$\rom;
$\sing(\XXf)=O\oplus\ssing$, where
$O$ is of type $\bX_{1,0}$,  $\bX_{1,1}$, or $\bZ_{11}$
and $\ssing$ is empty or $\bA_1$\rom;
\item\label{add.20}
if $\XXf$ has $20$ lines, then
%$\sing(\XXf)$ is either $\bX_{1,2}$ or as in item~\iref{add.18}.
either $\sing(\XXf)=\bX_{1,2}$ or $\sing(\XXf)$ is as in item~\iref{add.18}.
\endroster
In case~\iref{add.20}, the $20$ lines
consist of four coplanar
ones in the tangent plane $C_O\XXf$
and a generalized quadrangle $\GQ(3,1)$,
see \autoref{th.X} and \autoref{rem.X}.
\endaddendum

\autoref{add-non-K3-projective-double} is but one example
(\cf. also \autoref{thm-non-K3-projective})
of the
refined
statements that we can obtain using elementary algebraic topology and lattice
theory.
As yet another example, none of the possible $18$-line configurations
in \autoref{add-non-K3-projective-double}
is a subgraph of any of the $20$-line ones.

Even though
our main focus
is on the
normal non-K3 quartics,
for the sake of completeness
we consider non-isolated singularities as well (see
\autoref{S.nonisolated}).
These were extensively
studied in a number of classical papers, and our only contribution is the case
of a line of double points, which maximizes the number of lines.
An example of a maximal configuration was constructed in
\cite[Example 3.10]{gonzalezalonso}.

\theorem[see \autoref{proof-27-nonisolated}]\label{thm-27-nonisolated}
If
$\XXf \subset \PP^3(\CC)$ is a quartic surface
with non-isolated singularities that is not ruled by lines,
then $\XXf$ contains at most $27$ lines.
This bound is sharp, attained only at quartics with $\sing(\XXf)$ a line
of double points.
%Moreover, if  $\XXf$ attains the above bound then $\sing(\XXf)$ consists of a line of double points.
\endtheorem

Since
a quartic with an isolated triple point contains at most
$31$ lines (see \cite[p.~58]{rohn-dreifache-punkt}
and \autoref{th.T} for a slight refinement),
Theorems~\ref{thm-non-K3-projective-double} and \ref{thm-27-nonisolated}
imply the following immediate bound.

\corollary\label{thm-non-K3-projective}
Let $\XXf \subset \PP^3(\CC)$ be a non-$K3$ quartic surface
that is
not ruled by lines. Then $\XXf$ contains $31$ or at most $29$ lines.
Furthermore,
if
$\XXf$ contains more than $27$ lines, then
\roster*
\item
$\XXf$ is normal and has a single singular point~$O$, which is of type
$\bP_8:=\bT_{3,3,3}$
or \rom(for $27$ lines only\rom) $\bP_9:=\bT_{3,3,4}$\rom;
\item
$\XXf$ has $12$ pairwise distinct lines passing through~$O$.
\endroster
See \autoref{th.T}
 for a detailed description of the configurations of lines.
\endcorollary

The classification of complex normal projective quartic surfaces with at least one
non-simple singular point $O \in \XXf$  is found in \cite{degtyarev-90}.
By the results therein,  there exist 2523
%various
constellations of singularities on such surfaces, each of which contains at
most two non-simple points (with only a few configurations containing
two,
see \autoref{lem.rational}).
Pairs
$(\XXf, O)$, where $\XXf$ is a normal non-K3 quartics
and $O$ is a distinguished non-simple singular point,
 naturally split into four families (see \autoref{lem-coordinates-isolated}
and \autoref{conv.types}),
 depending on the type of the
point $O$ in Arnold's classification \cite{AVG1}.
 For most configurations of singularities, the minimal resolution $\tX$ of
$\XXf$ is rational
(see \autoref{lem.rational}), and
it is %{\5 this} case where
in this case that
$\XXf$ can contain many lines.

\table[t]
\caption{A summary
for normal rational quartics
(see \autoref{conv.table})}\label{tab.bounds}
\hrule height0pt
\[*
\def\ss{\vspace{2pt}}
\def\hr{\noalign{\hrule\ss}}
\def\minifirst#1{\hss$#1$\hss}
\def\={\text{-}}
\def\={}
\def\?{\rlap{\,?}}
\minitab\qquad\noalign{\hrule\ss}
O&\max&M&b_2(\tX)&E&\SSf_{\phantom0}&\text{bound}&\text{ref}\cr\noalign{\hrule\ss}
\bT\=&\ul{31}&\ul{31}&13&34&\bA_{11}&20\mapsto32&{\autoref{s.T}}
\cr
\bX\=&\ul{20}&\ul{20}&12&22&\bE_7\oplus\bA_3&\ul{16}\mapsto\ul{20}&{\autoref{s.X}}
\cr
\bJ\rlap{$\*$}\=&\ul{12}&27&11&14\rlap{$^*$}&\bE_8\oplus\bD_1&13\mapsto14&{\autoref{s.J*}, \autoref{sec-sharp-jstar}}
\cr
\bJ\=&?&48&11&16&\bD_9&16\mapsto16&{\autoref{s.J}}
\cr
\noalign{\hrule\ss}
%2\bX_{9}&20&12&2\bA_1\oplus\bD_1&\phantom0\ul{3}\mapsto11%&\checkmark
%\cr
%%\noalign{\hrule\ss}
%\bJ_{4,0}\*&27&2&\bD_1&\ul{1}\mapsto2%&\checkmark
%\cr
%%\noalign{\hrule\ss}
%2\bJ_{10}\*&27&3&\bD_1&\ul{1}\mapsto3%&\checkmark
%\cr
%\noalign{\hrule\ss}
\endminitab
\]
\endtable
The relevant data, both old and new, pertaining to rational non-$K3$ quartics
are collected in \autoref{tab.bounds}, which is explained
below.

\convention\label{conv.table}
The columns in \autoref{tab.bounds} are as follows:
\roster*
\item
$O$ is a reference to one of the four classes of surfaces under consideration, see
\autoref{conv.types} for the precise meaning of $\bT$, $\bX$ and $\bJ$;
\item ``$\max$''
is the maximal number of lines on a quartic $\XXf$ in
the
given class;
here
and elsewhere, the bounds known to be sharp are underlined;
\item
$M$ is the bound on the number of lines found in \cite{gonzalezalonso}
or references therein;
\item
$b_2(\tX)$ is the second Betti number of the minimal resolution of
singularities $\tX$ of a \emph{rational} representative of the
family, see \autoref{lem-betti-final};
\item
$E$ is the bound derived from N.~Elkies~\cite{elkies}, see
\autoref{cor.Elkies}; the $^*$ for the $\bJ\*$-series indicates the fact that
an extra trick has been used in the proof to reduce $17$ down to~$14$;
\item
$\SSf$ is the reduced intersection lattice of a rational representative, see
\autoref{s.Sigma};
throughout the paper, we let $\bD_1:=[-4]$;
\item
``bound'' is the lattice theoretic bound, in the form
(see the end of \autoref{s.Sigma})
\[*
\max\#\{\text{vectors in $\SSf$}\}\mapsto\max\#\{\text{lines}\}
\]
\item
``ref'' is a reference to the proofs and various refinements of the bounds.
\endroster
\endconvention

The lattice theoretic bounds in \autoref{tab.bounds} are
sometimes based on computer aided arguments carried out with the help of
\GAP~\cite{GAP4.13}.
(Though, most cases reduce to an analysis of Dynkin diagrams, which could
still be done manually.)
We list the bounds
based on
Elkies' brilliant idea~\cite{elkies} to emphasize the fact that
\emph{both \autoref{thm-non-K3-projective-double} and
the sharp bound	of at most $31$ lines in
\autoref{thm-non-K3-projective} can be shown without
any computer-aided
arguments} (see also \autoref{rem-elkies-rules}).

In this paper we  focus on
the \ser5- and \ser6-series, because
the sharp bound for the \ser3- (resp.\ \ser4-) case
%can be
are found in \cite{rohn-dreifache-punkt} (resp.\
\cite[Proposition~3.2]{gonzalezalonso}).
However,
we study the other two families as well, for we prove
various facts of the geometry of members of all families:
we compute the intersection lattices in all cases,
find all configurations allowed by the lattice theoretic constraints,
types of singularities \etc.\ (see, \eg, Theorems~\ref{th.T}, \ref{th.X} and  \autoref{rem.X}).
We also study the cases where
the resolution $\tX$ is no longer rational (see \autoref{tab.irrational}
in \autoref{lem.rational}).

There are several reasons to study the geometry of line configurations on
non-K3 quartics. Firstly, we
give precise answers to questions on a class of
varieties that  have been
 subject of great interest
for almost two hundred years.
Secondly, quartic surfaces and curves on them remain a useful tool to construct
various examples and test conjectures. Finally, we analyze several techniques introduced
in the last decade
to study configurations of low-degree curves on surfaces. Strangely enough,
 the bounds based on linear algebra \cite{elkies},
which we get almost for free, are just a few units worse than the
known sharp ones, whereas, quite unexpectedly, the lattice theoretic
bounds
turn out sharp in
many
cases.
It is worth emphasizing that,
in contrast to the $K3$-case,
we have no Torelli type theorem at our disposal,
so that the latter bounds boil down to
elementary algebraic
topology and lattice theory.

\subsection*{Contents of the paper}

In \autoref{S.geometry} we collect
basic facts on complex normal non-K3 quartics and discuss the bounds on
the number of lines that can be derived with Elkies' trick (see \autoref{s.Elkies}).
In particular, in \autoref{proof.main} we give the proof of \autoref{thm-non-K3-projective-double}.
Then, we discuss the lattice theoretic bounds for various classes of normal non-$K3$ quartics
(see \autoref{Sboundslattices}). Here,  algorithmic lattice theory combined with
the power of computer aided computations yields various insights into the geometry
and combinatorics of
line configurations on quartic surfaces.
\autoref{sec-sharp-jstar} contains a proof of the sharp bound for $\ser5$ (see \autoref{th.J*sharp});
it is based on \autoref{th.J*}.  Finally, in \autoref{S.nonisolated} we study
line configurations on non-normal quartics
and complete the
proof of \autoref{thm-27-nonisolated}.

\convention\label{conv.C}
We work over the field of complex numbers $\CC$;
therefore, from now on we abbreviate $\Cp3:=\Cp3(\C)$.
Throughout this paper, root lattices
	are assumed to be negative-definite,
	and rational curves  are assumed to be irreducible.
\endconvention

\subsection*{Acknowledgements}

We would like to thank Matthias Sch\"{u}tt for helpful discussions.

%%%%%%%%%%%%%%%%%%%%%%%%%%%%%%%%%%%%%%%%%%%%%%%%%%%%%%%%%%%%%%%%%%%%%%%
%%%%%%%%%%%%%%%%%%%%%%%%%%%%%%%%%%%%%%%%%%%%%%%%%%%%%%%%%%%%%%%%%%%%%%%
%%%%%%%%%%%%%%%%%%%%%%%%%  sec 2 geometry of non-k3 quartics
%%%%%%%%%%%%%%%%%%%%%%%%%%%%%%%%%%%%%%%%%%%%%%%%%%%%%%%%%%%%%%%%%%%%%%%%
%%%%%%%%%%%%%%%%%%%%%%%%%%%%%%%%%%%%%%%%%%%%%%%%%%%%%%%%%%%%%%%%%%%%%%%%

%\section{Geometry of non-$K3$ quartics with isolated singularities.}\label{S.geometry}
\section{Geometry of normal non-$K3$ quartics}\label{S.geometry}

In this section we collect various useful facts, mainly from
\cite{degtyarev-90,gonzalezalonso,ctc-wall-99}. In \autoref{proof.main}, we
prove \autoref{thm-non-K3-projective-double}.

\subsection{The classification}\label{s.classification}
The following version of the classification theorem
will play a crucial r\^{o}le in the sequel. A complete proof of the version
given
below can be found as the proof of
\cite[Lemma~6.3]{gonzalezalonso}, but it is inspired/based on the
considerations in  \cite{degtyarev-90} and the  unpublished preprint \cite{ctc-wall-99}.

\allowdisplaybreaks

\begin{lem}[see \cite{degtyarev-90}, {\cite[Thm~8.1]{ctc-wall-99}},  {\cite[Lemma~6.3]{gonzalezalonso}}] \label{lem-coordinates-isolated}
%	{\bf \rm \cite[Thm~8.1]{ctc-wall-99}}
Let $O \in \XXf$ be an isolated non-rational \emph{double}
point of an irreducible  quartic surface $\XXf \subset \PP^3(\CC)$.
	Then, there exist  polynomials
	$$
	Q_2=\sum_{i+j=2} q_{ij}x^iy^j \quad and \quad H_4=\sum_{i+j+k=4} h_{ijk}x^iy^jz^k
	$$
	such that,
in appropriate coordinates $\left(x:y:z:w\right)$ with
$O=\left(0:0:0:1\right)$,
the quartic  $\XXf$  is given
	by one of the following equations\rom:
	\begin{eqnarray}
		& & w^2z^2+wz Q_2 + H_4=0, \label{eq-q4} \\
		& & w^2z^2+w (y^3+zQ_2)+H_4 = 0  \label{eq-q56}
	\end{eqnarray}
%	$O=\left(x:y:z:w\right)=\left(0:0:0:1\right)$
and,
	if $\XXf$ is given by  \eqref{eq-q56}, then one of the following
holds\rom:	
\begin{eqnarray}
		& & h_{400} = h_{310} = q_{20} =0 \quad \text{and} \quad h_{301} \neq 0,\quad\text{or} \label{eq-q5} \\
		& & h_{400}=\frac{1}{4}q_{20}^2, \quad h_{310}=\frac{1}{2}q_{20}q_{11},  \quad q_{20} \neq 0,  \quad \text{and} \quad h_{301} = 0.  \label{eq-q6}
	\end{eqnarray}
	Moreover, if $\LLL \subset \XXf$ is a line that is not contained in the tangent cone $C_{O}\XXf$, then the coordinate change that leads to \eqref{eq-q56} can be chosen  in such a way that $\LLL$ is contained in the plane $\left\{x=0\right\}$.
\end{lem}

\convention\label{conv.types}
In the bulk of the paper, we consider pairs $(\XXf,O)$, where
\[
\aligned
&\text{$\XXf$ is a quartic with isolated singularities not ruled by lines and}\\
&\text{$O$ is a distinguished non-simple singular point of~$\XXf$}.
\endaligned
\label{eq.pair}
\]
(Equivalently, a quartic $\XXf$ with isolated singularities is not ruled by lines if
and only if it has no fourfold singular points, \ie, $\XXf$ is not a cone.)
We subdivide such pairs into the following four families:
\roster*
\item
the \ser3-series \cite[Theorem~4.6]{degtyarev-90},
if $O$ is a triple point;
\item
the \ser4-series \cite[Theorem~2.11]{degtyarev-90}, or \qser4 in \cite{gonzalezalonso,ctc-wall-99},
%if $\XXf$ is
given by \eqref{eq-q4};
\item
the \ser5-series \cite[Theorem~1.9]{degtyarev-90}, or \qser5 in
\cite{gonzalezalonso,ctc-wall-99},
%if $\XXf$ is
given by \eqref{eq-q56}, \eqref{eq-q5};
\item
the \ser6-series \cite[Theorem~1.7]{degtyarev-90}, or \qser6 in
\cite{gonzalezalonso,ctc-wall-99},
%if $\XXf$ is
given by \eqref{eq-q56}, \eqref{eq-q6}.
\endroster
We label the families above
according to the type of~$O$ in a very general member, which is
$\bP_8:=\bT_{3,3,3}$, $\bX_9:=\bX_{1,0}$, or $\bJ_{10}:=\bJ_{2,0}$ in the
notation of \cite{AVG1}.
%In other words,
We write, \eg, $\XXf\in\ser3$ to indicate that
$\XXf$ is in the \ser3-series.
%Finally, the
The difference between \ser5 and \ser6
(\emph{special} \vs. \emph{non-special} in
\cite{degtyarev-90}), both adjacent to $\bJ_{10}$,
is explained in
%\autoref{lem.one.line} below.
the next lemma.
\endconvention
\lemma[see {\cite[Lemma~2.6]{gonzalezalonso}}]\label{lem.one.line}
If $\XXf\in\ser5$ \rom(resp.\ $\XXf\in\ser6$\rom), then there is exactly one
\rom(resp.\ none\rom) line $\LLL_0\subset\XXf$ passing through~$O$.
Furthermore, for $\XXf\in\ser5$,
\roster
\item\label{i.<=1}
there is at most one other line $\LLX\subset\XXf$ intersecting~$\LLL_0$\rom;
\item\label{i.=1}
if present, $\LLX\cup\LLL_0$ constitutes the intersection
$C_{O}\XXf\cap\XXf$ \rom(which otherwise splits into $\LLL_0$ and a conic\rom)\rom;
\item\label{i.all.others}
as a consequence, if $\LLX$
%as in item~\rom{\iref{i.<=1}}
is present,
it intersects all other lines on~$\XXf$.
%\done
\endroster
Finally,
a line $\LLX$ as in item~\rom{\iref{i.<=1}} can be present only if
the
%non-simple
singular point $O$ of~$\XXf$ is of type $\bJ_{2,p}$, $p\ge0$.
\endlemma

\proof
Statements~\iref{i.<=1}--\iref{i.all.others} are proved in
\cite[Lemma~2.6]{gonzalezalonso}, where it is observed also that,
assuming~\eqref{eq-q56}, \eqref{eq-q5},
\[
\text{a line $\LLX$ as in \autoref{lem.one.line}\iref{i.<=1} is present if and only if $h_{220}=0$}.
\label{eq.LLX}
\]
For the last statement, we recall that, by \cite[Theorem 1.9]{degtyarev-90},
the point~$O$ of~$S$
is stably homeomorphic to the singular point $(1:0:0)$ of the
discriminant~$D$ of~$F$ with respect to~$w$. The latter is of type
\emph{other} than $\bJ_{2,p}$ if and only if the principal part
\[*
- 4h_{301}x^3z^3
 + (q_{11}^2 - 4h_{220})x^2z^2y^2
 + 2q_{11}xzy^4
 + y^6
\]
of~$D$ is a perfect cube. Comparing this to
\[*
\left(\frac23q_{11}xz + y^2\right)^3
 = \frac{8}{27}q_{11}^3x^3z^3
 + \frac43q_{11}^2x^2z^2y^2
 + 2q_{11}xzy^4
 + y^6,
\]
we find $2q_{11}^3=-27h_{301}\ne0$, see~\eqref{eq-q5}, and, hence,
$h_{220}=-q_{11}^2/12\ne0$.
\endproof

Most normal non-$K3$ quartics are rational.
If
we need to emphasize the fact that $\XXf$ is assumed rational, we write
$\XXf\in\ser4\rat$ or~$\ser5\rat$. The following lemma gives a characterization of such surfaces.
Here and below, $\mu$ stands for the (total) Milnor number of
(a set of) singularities.

\lemma[see, \eg, \cite{degtyarev-90,umezu81,umezu84}]\label{lem.rational}
A quartic~$\XXf$ as in \eqref{eq.pair} is rational
if and only if $O$ is the
only non-simple singular point of~$\XXf$ and its type is other than
$\bJ_{4,0}$ or~$\bX_{2,0}$. Otherwise, \ie, if $\sing(\XXf)$ is
\[*
\bX_{2,0}\oplus\ssing,\ \mu(\ssing)\le1,\qquad
2\bX_{9}\oplus\ssing,\ \mu(\ssing)\le3,\qquad
\bJ_{4,0},\qquad\text{\rm or}\qquad
2\bJ_{10}
\]
\rom(with $S\in\ser5$ in the last two cases\rom), then
$\XXf$ is
of elliptic ruled type \rom(\ie,  birationally equivalent to an elliptic ruled
surface\rom),
see
\autoref{tab.irrational}, where we still use \autoref{conv.lattice}, with the
$M$-column removed and the
$\max$-column replaced with the
\emph{exact} number of lines $\ls|\Fn\XXf|$.
\done
\table
\caption{Irrational quartics (see \autoref{lem.rational} and \autoref{conv.table})}\label{tab.irrational}
\[*
\def\ss{\vspace{2pt}}
\def\hr{\noalign{\hrule\ss}}
\def\minifirst#1{\hss$#1$\hss}
\def\={\text{-}}
\def\={}
\def\?{\rlap{\,?}}
\minitab\qquad\noalign{\hrule\ss}
\sing(\XXf)&\ls|\Fn\XXf|&b_2(\tX)&E&\SSf_{\phantom0}&\text{bound}
%&\text{ref}
\cr\noalign{\hrule\ss}
\bX_{2,0}\oplus\ssing&{4-\mu(\ssing)}&6&6&{\bD_4}&{\ul0\mapsto\ul4}
%&{\autoref{s.T}}
\cr
2\bX_{9}\oplus\ssing&{8-2\mu(\ssing)}&6&11&{\bD_4}&{\ul0\mapsto\ul8}
%&{\autoref{s.X}}
\cr
\bJ_{4,0}\=&\ul1&4&2&\bA_1\oplus\bD_1&{\ul1\mapsto2}
%&{\autoref{s.J*}}
\cr
2\bJ_{10}&{\text{$\ul2$ or $\ul3$}}&4&\ul3&\bA_1\oplus\bD_1&{\ul1\mapsto\ul3}
%&{\autoref{s.J}}
\cr
\noalign{\hrule\ss}
\endminitab
\]
\endtable
\endlemma

Finally, we recall the bounds on the number of lines that run through an isolated non-simple singularity on a quartic surface.

\lemma[see \cite{gonzalezalonso}]\label{lem.K-lines}
The maximal number of lines \emph{through the fixed non-simple singular
point~$O$} is
\[*
12\ \text{\rom(if $\XXf\in\ser3$\rom)},\quad
4\ \text{\rom(if $\XXf\in\ser4$\rom)},\quad
1\ \text{\rom(if $\XXf\in\ser5$\rom)},\quad\text{or}\quad
0\ \text{\rom(if $\XXf\in\ser6$\rom)}.
\]
\rom(If $\XXf\in\ser4$, these lines constitute the intersection
$C_O\XXf\cap\XXf$\rom; for $\XXf\in\ser5$, see \autoref{lem.one.line}.\rom)
Furthermore,
a very general representative of each family has exactly the
specified number of lines through~$O$.
\done
\endlemma

%%%%%%%%%%%%%%%%%%%%%%%%%%%%%%%%%%%%%%%%%%%%%%%%%%%%%%%%%%%%%%%%%
%%%%%%%%%%%%%%% 2.2 the minimal resolution of singularities
%%%%%%%%%%%%%%%%%%%%%%%%%%%%%%%%%%%%%%%%%%%%%%%%%%%%%%%%%%%%%%%%	

\subsection{The minimal resolution of singularities}\label{s.resolution}

We fix the notation $\sing(\XXf) = \{O, P_1, \dots P_r\}$,
where $O$ is the distinguished non-simple singular point;
occasionally we let $P_0:=O$.
In the sequel the minimal resolution of the quartic $\XXf$
(in the sense of \cite[p.~106]{bpv}) is denoted
as
\[*
\pi\: \tX \to\XXf(O)\to{\XXf},
\]
where $\XXf(O)$ is the
normalization of the proper transform of~$\XXf$ under the blow-up of
$\Cp3$ at~$O$,
and we use the shorthand
\[*
\KK:=K_{\tX},\qquad
\hh:=\pi^*[\Cp2],\qquad
E_i:=E(P_i):=\pi^{-1}(P_i),\quad i=0,\ldots,r,
\]
for the canonical class and hyperplane section of~$\tX$ and the
%and we put  $E_0 := \pi^{-1}(O)$ and  $E_i := \pi^{-1}(P_i)$, for $i=1, \ldots r$  to denote the
exceptional divisors of $\pi$. 	
Moreover,
by abuse of notation, we use $\sLLL$ to denote the proper transform under $\pi$ of
a  line $\sLLL \subset \XXf$.

\convention\label{conf.intersection}
Unless
explicitly stated otherwise, we say
that two lines \emph{intersect} if they do so in~$\tX$, not in $\XXf$;
in other words, we analyze the intersections after the singularities have been
resolved. (In particular, this applies to the dual adjacency graph
$\Fn\XXf$
introduced below.)
For example, all lines passing through~$O$ are considered
pairwise disjoint.
More generally, the dot $\cdot$ \emph{always} stands for
the intersection index in~$\tX$.
\endconvention

It is crucial that,
since $\deg\XXf=4$, the canonical class
$\KK$ is supported over the non-simple singular points of $\XXf$.
In fact, $-\KK$ is the
fundamental cycle of $E_0$ (or the sum thereof over all non-simple singular
points); in particular, this class is effective.

\begin{lem} \label{lem-minustwo}
If $\sLLL \subset \XXf$ is a line that does not \rom(resp.\ does\rom)
run through a non-simple singular point,
then
$\sLLL\cdot\KK=0$ and $\sLLL^2 = -2$ \rom(resp.\
$\sLLL\cdot\KK=\sLLL^2=-1$\rom) in $\tX$.
Moreover,
any class $\ell\in H_2(\tX)$ with
$
\sLLL\cdot\KK=\sLLL^2=-1 \mbox{ and } h \cdot \sLLL =1
$
is effective.
\end{lem}	

\proof
By the adjunction formula,
$\sLLL^2=-2-\sLLL\cdot\KK$.
As $-\KK$ is effective, we have $\sLLL\cdot\KK\le0$ and,
by the Riemann--Roch theorem, $\dim\ls|\ell|\ge-\ell\cdot\KK-1$. This implies the last
statement of the lemma and shows that $\XXf$ would be ruled by lines if it
had a line~$\ell$ with
$\ell\cdot\KK\le-2$.
\endproof

An important consequence taken for granted in the sequel is the fact that
lines have negative self-intersection. Hence, each class in $H_2(\tX)$ is
represented by at most one line and, instead of counting lines, we count
(or rather estimate the number of)
their homology classes.

For short, we will refer to a line $\sLLL\subset\XXf$ that does not
(resp.\ does) pass through a non-simple singular point of~$\XXf$ as a $(-2)$-line
\rom(resp.\ $(-1)$-line\rom). We use the notation
\begin{equation} \label{eq-fano-on-resolution}
\Fn\XXf=\FnZ\XXf\cup\FnK\XXf
\end{equation}
for the total dual adjacency graph of lines on~$\tX$ and its subgraphs of $(-2)$-
and $(-1)$-lines.

\remark\label{rem.vertical.plane}
The projection $\XXf\dashrightarrow\Cp2$ from~$O$ has degree~$1$ (if
$\XXf\in\ser3$) or~$2$ (otherwise). Hence, a plane through~$O$ may contain,
respectively, at most one or two $(-2)$-lines on~$\XXf$.
\endremark

The following statement is an immediate consequence of
the additivity of the topological Euler characteristic $\chtop$.
Recall that we put $P_0 := O$.

\lemma\label{lem-betti-prep}
Let $\XXf$ be a quartic with isolated singularities. Then
\[*
b_2(\tX)=22+4q(\XXf)+\sum\bigl(\chtop(E_i)-\mu(P_i)-1\bigr),
\]
the summation running over the \emph{non-simple} singular points
$P_i\in\sing(\XXf)$.
\done
\endlemma

\remark
In \autoref{lem-betti-prep} one can extend the summation to all points
$P_i\in\sing(\XXf)$ but the terms corresponding to the simple ones vanish.
It is this observation that constitutes the proof of the lemma.
\endremark

\corollary\label{lem-betti-final}
Let $(\XXf,O)$ be as in~\eqref{eq.pair}. If $\XXf$ is rational,
then $b_2(\tX)$ is as shown in \autoref{tab.bounds}.
Otherwise, it is as shown in \autoref{tab.irrational}.
\endcorollary

\proof
The exceptional divisors over most singularities involved look like
singular elliptic fibers except that they have self-intersection
\[
\KK^2=-3\ \text{(if $\XXf\in\ser3$)},\quad
-2\ \text{(if $\XXf\in\ser4$)},\quad \text{or}\quad
-1\ \text{(if $\XXf\in\ser5$ or \ser6)}.
\label{eq.K2}
\]
The minimal resolutions of corank~$2$ singularities are shown in
\autoref{tab.resolution}
(see also \autoref{rem.resolution}),
where we use both Kodaira's
notation and that in terms of affine Dynkin diagrams; for those of
corank~$3$ (the \ser3-series), we refer to~\cite{degt:quartics2}.

Arguing on the case-by-case basis, we can easily see that the difference
$\mu(O)-\chtop(E)\in\{8,9,10\}$
is constant within each series, provided that $O$ is
neither $\bX_{2,0}$ nor $\bJ_{4,0}$. (A more conceptual explanation of this
phenomenon is found in~\cite{degt:quartics2}, but it is difficult
to control the minimality of the resolution.)
Thus, for such points the statement follows
directly from \autoref{lem-betti-prep}.

\table[t]
\caption{Exceptional divisors of resolutions of elliptic singularities as
singular elliptic fibers (see \autoref{rem.resolution})}\label{tab.resolution}
\hbox to\hsize{\hss\mutab{$\bJ$, $\bE$: $\kappa^2=-1$, one $(-3)$-curve}
\bJ_{2,0}&\tA_0&\I_0&10\cr
\bJ_{2,1}&\tA_0^*&\I_1&11\cr
\bJ_{2,s}&\tA_{s-1}&\I_s&10+s\cr
\bJ_{3,s}&\tD_{s+4}&\I_s^*&16+s\cr
\bE_{12}&\tA_0^{**}&\II&12\cr
\bE_{13}&\tA_1^*&\III&13\cr
\bE_{14}&\tA_2^*&\IV&14\cr
\bE_{18}&\tE_6&\IV\*&18\cr
\bE_{19}&\tE_7&\III\*&19\cr
\bE_{20}&\tE_8&\II\*&20\cr
\endmutab
\quad
\mutab{$\bX$, $\bZ$: $\kappa^2=-2$, one $(-4)$-curve}
\bX_{1,0}&\tA_0&\I_0&9\cr
\bX_{1,1}&\tA_0^*&\I_1&10\cr
\bX_{1,s}&\tA_{s-1}&\I_s&9+s\cr
\bZ^1_{1,s}&\tD_{s+4}&\I_s^*&15+s\cr
\bZ^1_{11}&\tA_0^{**}&\II&11\cr
\bZ^1_{12}&\tA_1^*&\III&12\cr
\bZ^1_{13}&\tA_2^*&\IV&13\cr
\bZ^1_{17}&\tE_6&\IV\*&17\cr
\bZ^1_{18}&\tE_7&\III\*&18\cr
\bZ^1_{19}&\tE_8&\II\*&19\cr
\endmutab
\hss}
\bigskip
\hbox to\hsize{\hss\mutab{$\bY$, $\bW$: $\kappa^2=-2$, two $(-3)$-curves}
\bY^1_{r,s}&\tA_{r+s-1}&\I_{r+s}&9+r+s\cr
\bW_{1,s}&\tD_{s+4}&\I_s^*&15+s\cr
\bW^\sharp_{1,s}&\tD_{s+4}&\I_s^*&15+s\cr
\bW_{12}&\tA_1^*&\III&12\cr
\bW_{13}&\tA_2^*&\IV&13\cr
\bW_{17}&\tE_6&\IV\*&17\cr
\bW_{18}&\tE_7&\III\*&18\cr
\endmutab
\hss}
\endtable

The exceptional divisor over a point of type
$\bJ_{4,0}$ (resp.\ $\bX_{2,0}$)
splits into a smooth elliptic
curve~$E$
and a smooth rational curve~$R_1$
(resp.\ two smooth rational curves $R_1$, $R_2$),
so that
\[
E^2=-1\ \text{(resp.\ $-2$)},\qquad
R_i^2=-2,\qquad
E\cdot R_i=1,\qquad
R_1\cdot R_2=0.
\label{eq.J4}
\]
Thus, $\mu(O)-\chtop(E)=21$ (resp.\ $19$) in this case.
\endproof

\remark\label{rem.resolution}
In \autoref{tab.resolution},
%we leave out the self-intersections of the
%components, which are irrelevant for our problem.
%For the interested reader, we state
%the ultimate result below. It is these self-intersections that distinguish
%the types.
%the fundamental cycle~$\kappa$ has self-intersection $\kappa^2=-1$ for the
%$\bJ$- and $\bE$-series or $\kappa^2=-2$ for $\bX$, $\bY$, $\bZ$, $\bW$.
$\kappa^2$ is the self-intersection of the fundamental cycle.
If the exceptional divisor is irreducible, its self-intersection is
$\kappa^2$. Otherwise, in the first two tables, all but one vertices
represent rational $(-2)$-curves and one distinguished
\emph{simple} (\ie, one with coefficient~$1$ in~$\kappa$)
vertex represents a
rational $(\kappa^2-2)$-curve.
In the last table, two of the simple vertices represent $(-3)$-curves:
for~$\bY^1_{r,s}$, $\bW_{1,s}$, and \smash{$\bW^\sharp_{1,s}$},
they are at a distance of,
respectively, $r$, $2$, or $s+2$
from each other.

This computation gives us a description of the canonical divisor
$\KK=-\kappa$.
\endremark

\remark\label{rem.T.projection}
If $\XXf\in\ser3$, comparing \cite{degtyarev-90}
and~\cite{degt:quartics2}, one can easily see that the degree~$1$ projection
$\tX\to\Cp2$ contracts exactly $12$ rational curves: whenever $r$ of the $12$
lines through~$O$ (see \autoref{lem.K-lines})
collide to an $r$-fold line, an $\bA_{r-1}$ type singularity
of
%the blow-up
$\XXf(O)$
appears on (the proper transform of) that line. These $\bA$-type points are
all singularities of $\XXf(O)$,
\cf. \autoref{rem.K-lines} and \autoref{s.T} below.
\endremark

%%%%%%%%%%%%%%%%%%%%%%%%%%%%%%%%%%%%%%%%%%%%%%%%%%%%%%%%%%%%%%%%%%%
%%%%%%%%%%%% 2.3 elkies bound
%%%%%%%%%%%%%%%%%%%%%%%%%%%%%%%%%%%%%%%%%%%%%%%%%%%%%%%%%%%%%%%%%%%%

\subsection{Elkies' bound on the number of lines}\label{s.Elkies}

We are ready to state a simple
%, though still sufficient for \autoref{thm-non-K3-projective},
bound on the number of lines arising from
N.~Elkies~\cite{elkies}; in view of \cite{gonzalezalonso}, this bound
suffices to prove \autoref{thm-non-K3-projective-double}.

\corollary\label{cor.Elkies}
For a quartic~$\XXf$ as in \eqref{eq.pair} one has $\ls|\Fn\XXf|\le E$, where
$E$ is as given in Tables~\ref{tab.bounds}, \ref{tab.irrational}.
\endcorollary

\proof
We have $\ls|\Fn\XXf|\le\max\ls|\FnZ\XXf|+\max\ls|\FnK\XXf|$, and the second
term is bounded by \autoref{lem.K-lines} (which is to be doubled if there are
two non-simple points).
For the first term, project the
$(-2)$-lines to the vector space
\[*
V:=(\Q \hh\oplus\Q\KK)^\perp\subset H_2(\tX;\Q),\qquad
n:=\dim V=b_2(\tX)-2\le11.
\]
Since $K_{\tX}^2<0$ and $H_2(\tX)=\NS(\tX)$ is hyperbolic
by the Hodge index theorem,
$V$ is negative definite. The projection is
\[*
\LLL\mapsto l:=\LLL-\frac14\hh;
\]
the images of distinct lines are distinct and one has
$l^2=\qZ:=-9/4$ and
$l_1\cdot l_2=\qZ+2=1/4$ or $\qZ+3=3/4$ for $l_1\ne l_2$.

Rescale the form on~$V$ by $-4/9$,
 so that $V$ be positive definite, $l^2=1$,
and the products $l_1\cdot l_2$ take but two values $\tau_1=1/9$,
$\tau_2=-1/3$. Since $\tau_1+\tau_2\le0$ and $1+n\tau_1\tau_2>0$,
from~\cite{elkies} we find that
\[*
\ls|\FnZ\XXf|\le\frac{(1-\tau_1)(1-\tau_2)n}{1+\tau_1\tau_2n}=\frac{32n}{27-n}.
\]
Together with \autoref{lem-betti-final}
this concludes the proof for the rational surfaces in all series
except $\bJ\*$, for which we obtain
$\ls|\Fn\XXf|\le17$. To improve this last bound, we ignore
(if present) the only line~$\LLX$
intersecting~$\LLL_0$, see~\autoref{lem.one.line}, and project the rest to the
smaller space $\LLL_0^\perp\subset V$;
this time $\qZ=-2$ and
$\tau_1=0$, $\tau_2=-1/2$.
Upon applying~\cite{elkies}, we add $1$ to the result to
account for~$\LLX$.

For the few irrational surfaces~$\XXf$ (see \autoref{lem.rational}),
we use the smaller values of
$b_2(\tX)$ given by \autoref{tab.irrational} and change~$V$ to the orthogonal
complement of~$h$ and
the subspace generated by
\emph{all components} of the exceptional divisors over
the non-simple points, see~\eqref{eq.J4}.
\endproof

%\begin{lem}[{see \cite[Proposition~3.2]{gonzalezalonso}}]\label{lem-q4}
%	If $\XXf$ belongs to \ser4 then it contains at most $20$ lines.\mnote{this
%should go elsewhere, to the proof of the main theorem}
%\done
%\end{lem}
%%\begin{proof} See  \cite{gonzalezalonso}.
%%\end{proof}

%\proof[Proof of \autoref{thm-non-K3-projective}]
%\endproof

\subsection{Proof of \autoref{thm-non-K3-projective-double}}\label{proof.main}

By \autoref{cor.Elkies}, both for $\XXf \in \ser5$ and for $\XXf \in \ser6$
we have  $\ls|\Fn\XXf|\le 16$. On the other hand, for $\XXf \in \ser4$ we
have the sharp bound $\ls|\Fn\XXf|\le 20$ by \cite[Proposition~3.2]{gonzalezalonso} (see $M$ in
\autoref{tab.bounds}). Finally, by  \autoref{cor.Elkies}, the maximal value of $20$ lines is never attained when $\XXf$ has two non-simple points
(see $E$ in
\autoref{tab.irrational}).
This completes the proof.
\qed

\section{Bounds on the number of lines \via\ lattices}\label{Sboundslattices}

In the bulk of this section, $(\XXf,O)$ is a \emph{rational} pair as
in~\eqref{eq.pair}; in particular, by \autoref{lem.rational},
$O$ is the only non-simple singular point.
An exception is \autoref{s.irrational}.

\subsection{The reduced intersection lattice}\label{s.Sigma}

According to \autoref{lem-betti-final}, the second Betti number
$b_2(\tX)$ stays constant over the rational surfaces within each of the
four families in \autoref{conv.types}.
Hence, so does the intersection lattice
$H_2(\tX):=H_2(\tX;\Z)$ and the pair
of classes $\KK,\hh\in H_2(\tX)$. We define the
\emph{reduced intersection lattice} of~$\tX$ as
\[*
\SSf:=(\Z \hh\oplus\Z\KK)^\perp\subset H_2(\tX);
\]
it is a negative definite lattice of rank $b_2(\tX)-2$.

\lemma\label{lem.Sigma}
The lattice $\SSf$ is as given in \autoref{tab.bounds}.
\endlemma

\proof
The intersection lattice $H_2(\tX)$ of the rational surface~$\tX$ is
$[1]\oplus\bH_{b_2(\tX)-1}$, where we fix the notation
\[*
\bH_n=n\bH_1:=\bigoplus_{i=1}^n\Z e_i,\quad e_i^2=-1,
\]
for the standard negative definite Euclidean lattice. We have $\hh^2=4$,
$\hh\cdot\KK=0$, and $\KK^2$ is given by~\eqref{eq.K2}. Since also
$K_X=w_2(X)\bmod2$ for any algebraic surface~$X$, we conclude that $\SSf$ is
an \emph{even} negative definite lattice of rank $b_2(\tX)-2$. Its \emph{genus}
(equivalently, discriminant $\discr\SSf$) is easily computed using
Nikulin~\cite{nikulin:forms}, and it is indeed as in the statement.

To find all representatives of each genus, we use~\cite{nikulin:forms} again
and show that $\SSf=T^\perp\subset L$ for an appropriate
\emph{characteristic} sublattice~$T$ of an odd negative definite unimodular
lattice~$L$ of rank~$12$. There are three such lattices, \viz.
\[*
L=\bH_{12},\quad
\bE_8\oplus\bH_4,\quad\text{or}\quad
\bD_{12}^+,
\]
where the latter is an index~$2$ extension of $\bD_{12}$
\emph{other than $\bH_{12}$}.
Thus, it remains to indicate~$T$, list all representatives of each genus, and
use geometric insight to select the ``correct'' one. Essentially this is
done in~\cite[Theorem 4.2]{degt:quartics2}. Below, at the beginning of each of
\autoref{s.T}--\autoref{s.J}, we complete the proof of the lemma
by providing a simpler
geometric argument.
\endproof

\convention\label{conv.lattice}
Given an even negative definite lattice $\SSf$, a discriminant class
$\Ga\in\discr\SSf=\SSf\dual\!/\SSf$,
and a rational number $q=\Ga^2\bmod2\Z$, we
use the following notation:
\roster*
\item
$\Vec(\SSf,\Ga,q)$ is the set of $v\in\SSf\dual$ such that $v^2=q$ and
$v\bmod\SSf=\Ga$;
\item
$\bnd(\SSf,\Ga,q)$ is the maximal cardinality of a subset $V\subset\Vec(\SSf,\Ga,q)$
such that
\[
u\cdot v\in\{q+2,q+3\}\quad\text{for any pair $u\ne v$ in $V$}.
\label{eq.V}
\]
\endroster
We abbreviate
\[*
\VecK(\SSf,\Ga),\ \bndK(\SSf,\Ga)\quad\text{or}\quad
\VecZ(\SSf,\Ga),\ \bndZ(\SSf,\Ga)
\]
if $-2<q\le0$ or $-4<q\le-2$, respectively, which are the two relevant cases.
\endconvention

Consider the orthogonal projection $p\:H_2(\tX)\to\SSf\dual$. The images of
lines are
\[*
\ell\mapsto l:=\ell-\frac14\hh\quad\text{for a $(-2)$-line},\qquad
\ell\mapsto l:=\ell-\frac14\hh+\kappa\KK\quad\text{for a $(-1)$-line},
\]
where $\kappa:=1/\KK^2$, see~\eqref{eq.K2}. The following
statement is immediate.

\lemma\label{lem.proj}
The projection $p\:H_2(\tX)\to\SSf\dual$ has the following properties\rom:
\roster*
\item
the images of distinct lines are distinct\rom;
\item
the images $l$ of all $(-2)$-lines are in the same class
$\eta\in\discr\SSf$\rom; one has $l^2=\qZ:=-9/4$\rom;
\item
the images $l$ of all $(-1)$-lines are in the same class
$\lambda\in\discr\SSf$\rom; one has $l^2=\qK:=-\kappa-5/4$\rom;
\item
for the images $l_0$, $l_1$, $l_2$ of a $(-1)$-line $\LLL_0$ and $(-2)$-lines
$\LLL_1\ne\LLL_2$ one has
\[*
l_1\cdot l_2=\LLL_1\cdot\LLL_2-\frac14,%\in\{\qZ+2,\qZ+3\},
\qquad
l_1\cdot l_0=\LLL_1\cdot\LLL_0-\frac14;%\in\{\qZ+2,\qZ+3\}.
\]
thus, the intersections take values in $\{\qZ+2,\qZ+3\}$.
%\item\label{pp.FnZ}
%as a consequence, $p$ bijects $\FnZ\XXf$ onto a subset
%$V\subset\VecZ(\SSf,\eta)$ as in \autoref{conv.lattice}\rom;
%\item\label{pp.FnK}
%likewise, $p$ injects $\FnK\XXf$ into $\VecK(\SSf,\lambda)$.
\endroster
\endlemma

We reserve the notation $\eta,\lambda$ and $\qZ,\qK$ introduced
in \autoref{lem.proj} for the rest of this section. The classes
$\eta,\lambda\in\discr\SSf$ are such that
\[*
\eta^2=-1/4\bmod2\Z,\quad
\lambda^2=-\kappa-5/4\bmod2\Z,\quad
\eta\cdot\lambda=-1/4\bmod\Z.
\]
In each lattice $\SSf$ considered below,  a
pair of classes with these properties (a single class $\eta$ in case \ser6) is unique up to
$\OG(\SSf)$ (in fact, up to $\pm1$). Hence, we assume them known and fixed.

As follows from \autoref{lem.proj}, the projection establishes bijections
\[
\FnZ\XXf\overset\cong\longto\VZ(\XXf)\subset\VecZ(\SSf,\eta),\qquad
\FnK\XXf\overset\cong\longto\VK(\XXf)\subset\VecK(\SSf,\lambda).
\label{eq.sets}
\]
By \autoref{lem.proj} again, the set $\VZ(\XXf)$ satisfies \eqref{eq.V};
hence,
\[*
\ls|\FnZ\XXf|\le\bndZ(\SSf,\eta).
\]
It is this purely arithmetical
bound that is denoted as ``$\max\#\{\text{vectors in $\SSf$}\}$'' in
\autoref{conv.table} and used in Tables~\ref{tab.bounds},
\ref{tab.irrational}. The other integer in the column ``bound'' in \autoref{tab.bounds}
is obtained by merely adding
$\max\ls|\FnK\XXf|$, see \autoref{lem.K-lines}.

Besides, we have the sets
\[*
\Eset(\XXf)\subset\Eall(\XXf)\subset\VecZ(\SSf,0)
\]
of the \emph{exceptional $(-2)$-divisors} on~$\tX$,
\ie, smooth rational
$(-2)$-curves that are orthogonal to~$\hh$, and, respectively, all positive roots in the
root lattice generated by $\Eset(\XXf)$.
These are either the irreducible components of the exceptional divisors
$E_i$, $i\ge1$, over the simple singular points of~$\XXf$ or rational
components of $E_0$ \emph{orthogonal to $\KK$}. We call~$\XXf$
\emph{relatively smooth}
if $\Eset(\XXf)=\varnothing$.

Finally, we consider the set
\[*
\Cset(\XXf)\subset\Q\VZ(\XXf)^\perp\cap\SSf\dual\subset\SSf\dual
\]
of the images under~$p$ of the rational components of~$E_0$.
A component $C$ of square $C^2=-2$, $-3$, $-4$, or~$-5$ (\cf.
\autoref{rem.resolution}) projects to a vector $c\in\Cset(\XXf)$ of square
\[
c^2=C^2-\kappa(C^2+2)^2.
\label{eq.Cset}
\]
By the construction, $E_0$ is irreducible if and only if
$\Cset(\XXf)=\varnothing$.

%%%%%%%%%%%%%%%%%%%%%%%%%%%%%%%%%%%%%%%%%%%%%%%%%%%%%%%
%%%%%%%%%%%%%  3.2
%%%%%%%%%%%%%%%%%%%%%%%%%%%%%%%%%%%%%%%%%%%%%%%%%%%%%%%%%

\subsection{Geometric restrictions}\label{s.geometric}

Unlike the case of $K3$-quartics, the lattice
$H_2(\tX)=\NS(\tX)\ni\hh,\KK$
does not give us full control over the configuration of lines. Below, we
state a few simplest restrictions on the sets introduced that arise from the
geometry of quartics.

\lemma\label{lem.ex}
The set $\Eall(\XXf)$ has the following properties\rom:
\roster*
\item
$e\cdot l\ge0$ for each $e\in\Eall(\XXf)$ and $l\in\VZ(\XXf)\cup\VK(\XXf)$\rom;
\item
if $e\cdot l_1=e\cdot l_2=1$ for $e\in\Eall(\XXf)$, $l_1,l_2\in\VZ(\XXf)$,
then $l_1\cdot l_2=\qZ+2$.
\endroster
\endlemma

\proof
For the former, all divisors involved are effective and without common
components. For the latter, geometrically the corresponding lines $\LLL_1$,
$\LLL_2$ intersect \emph{in~$\XXf$} at the singular point of~$\XXf$ to which $e$ contracts;
hence, they have no other intersection points.
\endproof

\lemma[the triangle property]\label{lem.K3}
Given three classes $l_1,l_2,l_3\in\VZ(\XXf)$ that satisfy
$l_i\cdot l_j=\qZ+3$ for all $1\le i<j\le3$ \rom(in other words, the corresponding lines
$\LLL_1,\LLL_2,\LLL_3$
intersect\rom),
 either
\roster*
\item
there is a vector $e\in\Eall(\XXf)$ such that $e\cdot l_i=1$ for some
$i=1,2,3$, or
\item
there is a fourth class $l_4\in\VZ(\XXf)$ such that
$l_4\cdot l_i=\qZ+3$, $i=1,2,3$.
\endroster
\endlemma

\proof
The three lines as in the statement are coplanar;
let~$\Pi$ be the plane spanned by these lines, and let $\LLL_4$ be the fourth
component of the degree~$4$ curve $\XXf\cap\Pi$. By
\autoref{rem.vertical.plane}, $\Pi\not\ni O$ and $\LLL_4$ is a $(-2)$-line.

If $\LLL_4=\LLL_i$ for some $i=1,2,3$, then $\Pi$ is tangent to~$\XXf$
along~$\LLL$ and, hence, there is a singular point of~$\XXf$ on $\LLL=\LLL_i$ (as
otherwise the normal bundle of~$\LLL$ in~$\XXf$ would be that in~$\Pi$,
implying $\LLL^2=1$).
Otherwise, $\LLL_4$ intersects each of $\LLL_i$ \emph{in $\XXf$} and the
intersection points survive to~$\tX$ unless they are singular for~$\XXf$.
\endproof

\remark\label{rem.K-lines}
By \autoref{lem-minustwo},
every divisor $\LLL\in H_2(\tX)$ with the property that
$\ell\cdot\KK=\ell^2=-1$ and $\ell\cdot\hh=1$ is effective. Hence, each
vector $l\in\Vec(\SSf,\lambda)$ is the image of a unique ``line'' $\LLL$
on~$\XXf$ through~$O$. However, we cannot assert that this line $\LLL$ is
irreducible; it may happen that $\LLL=\LLL'+e$ for another line $\LLL'$
through~$O$ (possibly, still reducible) and an exceptional divisor
$e\in\Eset(\XXf)$ such that $e\cdot\LLL'=1$
(\cf. a
similar discussion of the relation between ``multiple'' lines in $C_O\XXf$ and
singular points in \cite[Lemmas~2.6 and~4.2]{degtyarev-90} and
in \autoref{rem.T.projection}).

\endremark

\remark\label{rem.brute.force}
In \autoref{s.T}--\autoref{s.J} below, we provide a
%detailed
combinatorial
description of the sets $\VZ(\XXf)$, $\VK(\XXf)$, and
$\Eset(\XXf)\subset\Eall(\XXf)$ which should suffice to derive most of our
classification statements manually. However, in most cases we choose to save
time/space and apply brute force, using \GAP~\cite{GAP4.13}.
Namely,
we
\roster
\item\label{brute.1}
list the subsets $V\subset\VecZ(\SSf,\eta)$ satisfying~\eqref{eq.V},
\item\label{brute.2}
for each~$V$, compute the \emph{maximal} subset $\Emax$ given by \autoref{lem.ex}, and
\item\label{brute.3}
use
this set $\Emax$ to eliminate the sets~$V$ violating \autoref{lem.K3}.
\endroster
Thus, we obtain a reasonably short list of \emph{candidates} for the
configuration of lines on a quartic. We never assert that all candidates are
realizable: the realizability is to be established by explicit examples.
\endremark

%%%%%%%%%%%%%%%%%%%%%%%%%%%%%%%%%%%%%%%%%%%%%%%%%%%%%
%%%%%%%%%%%%%  3.3 the t-series
%%%%%%%%%%%%%%%%%%%%%%%%%%%%%%%%%%%%%%%%%%%%%%%%%%%%%%

\subsection{The $\bT$-series}\label{s.T}

To complete the proof of \autoref{lem.Sigma}, we use $T=[-12]$, arriving at
the three candidates,
\[*
\SSf=\bA_{11},\quad
\bE_8\oplus\bA_2\oplus\bD_1,\quad\text{or}\quad
\bD_9\oplus\bA_2,
\]
with $\bigl|\Vec(\SSf,\Gl)\bigr|=12$, $3$, or $0$, respectively.
By \autoref{lem.K-lines} and \eqref{eq.sets},
only $\bA_{11}$ may serve as $\SSf$.

One can take for~$\lambda$ one of the two standard generators (those of
square $-11/12\bmod2\Z$) of the group $\discr\bA_{11}\cong\Z/12$; then,
$\eta=-3\lambda$. Consider the lattice $\bH_{12}$, denote by
$\IS:=\{1,\ldots,12\}$ the index set, and, for a subset $s\subset\IS$, let
$\vv{s}:=\ls|\IS|\1\sum_{i\in s}e_i\in\bH_{12}\otimes\Q$.
Then, $\bA_{11}$ is
$\vv\IS^\perp\subset\bH_{12}$, and we have
\[*
\aligned
\VecK(\bA_{11},\lambda)&=
 \bigl\{12\cdot\vv{\{p\}}-\vv\IS\bigm|p\in\IS\bigr\},\\
\VecZ(\bA_{11},\eta)&=
 \bigl\{3\cdot\vv\IS-12\cdot\vv{s}\bigm|s\subset\IS,\ \ls|s|=3\bigr\},\\
\VecZ(\bA_{11},0)&=
 \bigl\{e_i-e_j\bigm|i,j\in\IS,\ i\ne j\bigr\}.
\endaligned
\]
Thus, by~\eqref{eq.sets}, we can
%identify the $(-1)$- (resp.\ $(-2)$-) lines
%on~$\tX$ with points $p\in\IS$ (resp.\ $3$-element sets
%$s\subset\IS$). Therefore, we regard $\FnK\XXf=\VK(\XXf)$ as a subset of~$\IS$
%and $\FnZ\XXf=\VZ(\XXf)$, as a set of $3$-element subsets $s\subset\IS$.
\roster*
\item
identify the set $\FnK\XXf=\VK(\XXf)$ of $(-1)$-lines
with a subset of $\IS$, and
\item
identify the set $\FnZ\XXf=\VZ(\XXf)$ of $(-2)$-lines
with a certain collection of $3$-element subsets $s\subset\IS$.
\endroster
Under this identification, we have:
\roster*
\item
$\ls|s_1\cap s_2|\le1$ for any pair $s_1\ne s_2$ in~$\VZ(\XXf)$
as a restatement of~\eqref{eq.V};
\item
a $(-1)$-line $\pp\in\IS$ and $(-2)$-line $s\subset\IS$ intersect in~$\tX$
if and only if $\pp\in s$;
\item
two $(-2)$-lines $s_1,s_2\subset\IS$ intersect in $\tX$
if and only if $s_1\cap s_2=\varnothing$.
\endroster
According to \autoref{rem.T.projection} (\cf. also \autoref{rem.K-lines}), we
can re-index~$\IS$ so that
\roster*
\item
%the set
$\VK(\XXf)$ consists of $r\le12$ points
$\pp_1<\pp_2<\ldots<\pp_r=12$; let $\pp_0:=0$;
\item
%the set
$\Eset(\XXf)$ consists of the $(12-r)$ divisors $e_i-e_{i+1}$,
$i\in\IS\sminus\VK(\XXf)$;
\item
if a $(-2)$-line $s\in\VZ(\XXf)$ contains a point $\pp\in(\pp_{k-1},\pp_k]$,
it contains the interval $(\pp_{k-1},\pp]$.
\endroster
Geometrically, a
%$(-2)$-line
set $s\in\VZ(\XXf)$ has non-empty intersection with
$(\pp_{k-1},\pp_k]$ if and only if the
corresponding lines $s$ and $\pp_k$ intersect \emph{in $\XXf$}.
The following lemma is an immediate consequence of the above
description (\cf. \cite[$\S$8]{rohn-dreifache-punkt}).

\lemma\label{lem.T.smooth}
For a quartic $\XXf\in\ser3$ one has
\[*
\ls|\FnK\XXf|+\ls|\Eset(\XXf)|=12.
\]
In particular, $\XXf$ is relatively smooth if and only if
it has the maximal number of $(-1)$-lines, \ie,
$\ls|\FnK\XXf|=12$ \rom(see \autoref{lem.K-lines}\rom).
\done
\endlemma

For the \ser3-series, we merely reconfirm arithmetically
%slightly refine
the sharp upper bound
$\ls|\Fn\SSf|\le31$ found in \cite[p.~58]{rohn-dreifache-punkt}
and restate, in the modern terms, a few results of~\latin{loc.\ cit.}\
concerning the large configurations of lines.
%~\cite{gonzalezalonso}.

\theorem[\cf. \cite{rohn-dreifache-punkt,gonzalezalonso}]\label{th.T}
If $\XXf\in\ser3$, then
$\ls|\Fn\XXf|=31$ or $\ls|\Fn\XXf|\le29$. These bounds are sharp.
Furthermore, if $\ls|\Fn\XXf|\ge28$, then $O$ is the only singular point
of~$\XXf$, one has
$\ls|\FnK\XXf|=12$, and either
\roster*
\item
$\FnZ\XXf$ is a generalized quadrangle $\GQ(3,1)$ or a $1$- or $3$-vertex
extension thereof, see \autoref{fig.GQ}, and $O$ is of type $\bP_8$, or
\item
$\FnZ\XXf$ is one of $U'_{16}$, $U''_{16}$ in \autoref{fig.U} and $O$ is of type
$\bP_8$ or $\bP_9$.
\endroster
\endtheorem

\remark\label{rem.T.sharp}
Unlike
most other statements in this paper, in \autoref{th.T} we also assert
the existence of all configurations/singularities described. It should not be
difficult to modify the proof to obtain a full deformation classification
in the modern language.
We do not engage into this part,
partially because it is also discussed  in \cite{rohn-dreifache-punkt}.
\endremark

\remark\label{rem.T.figures}
In Figures~\ref{fig.GQ}, \ref{fig.U}, the $(-2)$ lines are the columns: each
$(-2)$-line~$\LLL$ is interpreted as a $3$-element subset of $\FnK\XXf$,
\viz. the set of the $(-1)$-lines that $\LLL$ intersects.
Two $(-2)$-lines
intersect each other if and only if the corresponding subsets are disjoint.
In \autoref{fig.GQ}, $V_{n}$, $n=16,17$, is made of the first $n$ columns.
\endremark

\proof[Proof of \autoref{th.T}]
Interpreting
the elements of $\VecZ(\bA_{11},\eta)$ as $3$-element subsets
$s\subset\IS$, by \cite[Lemma~4.2]{degt:sextics} we have
$\bndZ(\bA_{11},\eta)\le20$. By brute force we confirm that this bound is
sharp.
Next, we follow \autoref{rem.brute.force} (\cf. \autoref{rem.T.brute.force}
below)
and select the sets $V\subset\VecZ(\bA_{11},\eta)$ satisfying
\autoref{lem.K3}.
There are (with the subscript always indicating the
cardinality):
\roster*
\item
three sets $V_{16}\cong\GQ(3,1)\subset V_{17}\subset V_{19}$  in
\autoref{fig.GQ}, all with $\Emax=\varnothing$,
\item
two sets $U'_{16}$ and $U''_{16}$  in \autoref{fig.U},
also with $\Emax=\varnothing$, and
\item
three more sets $U'_{17}$, $U''_{17}$, and $W_{17}$, this time with
$\Emax\ne\varnothing$.
\endroster
The last three sets are ruled out by \autoref{lem.T.smooth}, which leaves
room for at most one exceptional $(-2)$-divisor. Trying $1$-element
subsets
$\Eset\subset\Emax$ one-by-one, in each case we find a contradiction to
\autoref{lem.K3}.
Thus, we conclude that $\XXf$ is relatively smooth and has $12$ pairwise
distinct $(-1)$-lines
and there are but five candidates, \viz. those in Figures~\ref{fig.GQ}
and~\ref{fig.U},
for the configuration $\FnZ\XXf$.

\figure
\[*
\tconfig
 . . . . . . 1 1 1 . . . . . . 1 . . . \cr
 1 . . 1 . . . . . . 1 . 1 . . . 1 . . \cr
 . 1 . . 1 . . . . . . 1 . 1 . . 1 . . \cr
 . . 1 . . 1 . . . 1 . . . . 1 . 1 . . \cr
 . . . . . . . . . 1 1 1 . . . 1 . . . \cr
 1 . . . . 1 1 . . . . . . 1 . . . 1 . \cr
 . 1 . 1 . . . 1 . . . . . . 1 . . 1 . \cr
 . . 1 . 1 . . . 1 . . . 1 . . . . 1 . \cr
 . . . . . . . . . . . . 1 1 1 1 . . . \cr
 1 . . . 1 . . 1 . 1 . . . . . . . . 1 \cr
 . 1 . . . 1 . . 1 . 1 . . . . . . . 1 \cr
 . . 1 1 . . 1 . . . . 1 . . . . . . 1 \cr
\endtconfig
\]
\caption{The sets $V_{16}\cong\GQ(3,1)\subset V_{17}\subset V_{19}$ (see \autoref{rem.T.figures})}\label{fig.GQ}
\endfigure

\figure
\[*
\tconfig
 1 . . . 1 . . . . . . 1 . 1 . . \cr
 . 1 . . . 1 . . . . . . 1 . 1 . \cr
 . . 1 . . . 1 . . . . . . 1 . 1 \cr
 . . . 1 1 . . 1 . . . . . . 1 . \cr
 1 . . . . 1 . . 1 . . . . . . 1 \cr
 . 1 . . 1 . 1 . . 1 . . . . . . \cr
 . . 1 . . 1 . 1 . . 1 . . . . . \cr
 . . . 1 . . 1 . 1 . . 1 . . . . \cr
 1 . . . . . . 1 . 1 . . 1 . . . \cr
 . 1 . . . . . . 1 . 1 . . 1 . . \cr
 . . 1 . . . . . . 1 . 1 . . 1 . \cr
 . . . 1 . . . . . . 1 . 1 . . 1 \cr
\endtconfig\qquad
\tconfig
 1 1 . . 1 . . 1 . . . . . . . . \cr
 1 . . 1 . 1 . . . . . . . . . 1 \cr
 1 . 1 . . . . . . 1 1 . . 1 . . \cr
 . 1 . . . 1 . . . 1 . . . . . . \cr
 . 1 . . . . 1 . 1 . 1 . . . . . \cr
 . . 1 . 1 . . . . . . . . . . 1 \cr
 . . . . . 1 . . 1 . . . . 1 1 . \cr
 . . . 1 . . 1 . . . . 1 . 1 . . \cr
 . . . . . . . 1 . 1 . 1 . . 1 1 \cr
 . . 1 . . . . . 1 . . 1 1 . . . \cr
 . . . 1 . . . 1 . . 1 . 1 . . . \cr
 . . . . 1 . 1 . . . . . 1 . 1 . \cr
\endtconfig
\]
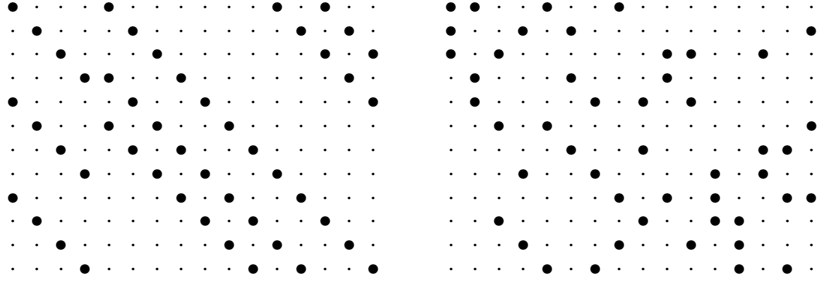
\caption{The sets $U'_{16}$ and $U''_{16}$ (see \autoref{rem.T.figures})}\label{fig.U}
\endfigure

The $(-1)$-lines define $12$ pairwise distinct smooth
points $p_1,\ldots,p_{12}$
(cut off by a plane quartic)
on the plane cubic~$E_0$, and the
$(-2)$-lines are collinearities $(p_i,p_j,p_k)$ of these
points.
In other words, $(-2)$-lines are
relations of the form $p_i+p_j+p_k=0$, $\{i,j,k\}\in V$,
in the group law on~$E_0$, which is
\roster
\item\label{group.A0}
$G:=(\R/\Z)^2$ if $E_0$ is smooth, \ie,
$O$ is of type~$\bP_8:=\bT_{3,3,3}$;
\item\label{group.A0*}
$G_m:=(\R/\Z)\times\R$ if $E_0$ is nodal, \ie,
$O$ is of type~$\bP_9:=\bT_{3,3,4}$;
\item\label{group.A0**}
$G_a:=\R^2$ if $E_0$ is cuspidal, \ie,
$O$ is of type~$\bQ_{10}$;
\item\label{group.A1}
$G_m\times(\Z/2)$ if $E_0$ is of type $\tA_1$, \ie,
$O$ is of type~$\bT_{3,4,4}$;
\item\label{group.A1*}
$G_a\times(\Z/2)$ if $E_0$ is of type $\tA_1^*$, \ie,
$O$ is of type~$\bS_{11}$;
\item\label{group.A2}
$G_m\times(\Z/3)$ if $E_0$ is of type $\tA_2$, \ie,
$O$ is of type~$\bT_{4,4,4}$;
\item\label{group.A2*}
$G_a\times(\Z/2)$ if $E_0$ is of type $\tA_2^*$, \ie,
$O$ is of type~$\bU_{12}$.
\endroster
Thus, each set~$V$ to be considered gives rise to a
%certain
``system of linear equations''
on some of these groups; we denote by~$M$ its ``matrix'',
which is essentially given
by the figures.
We emphasize that we are only interested in
\[
\text{solutions to the system with all $p_1,\ldots,p_{12}$ pairwise distinct}.
\label{eq.T.distinct}
\]
The existence of such a solution is also \emph{sufficient} for the
realizability of~$V$: in all cases considered, the extra relation
$p_1+\ldots+p_{12}=0$ (the fact that the points are cut off by a quartic
curve) follows from the system.

For $V=V_{19}$, $U'_{16}$, or $U''_{16}$, we have
\roster*
\item
$\Q V=\SSf\otimes\Q$, \ie, $E_0$ is irreducible, see~\eqref{eq.Cset},
leaving \iref{group.A0}--\iref{group.A0**} only;
\item
$\rank M=11$, \ie, there are no solutions in the torsion free group \iref{group.A0**}.
\endroster
For $U'_{16}$, $U''_{16}$, we do find solutions
satisfying~\eqref{eq.T.distinct} in $\R/\Z$; hence, also in~\iref{group.A0}
and \iref{group.A0*}. For~$V_{19}$, all invariant factors of~$M$ divide~$6$;
hence, any solution in~\iref{group.A0*} is $6$-periodic and cannot
satisfy~\eqref{eq.T.distinct}. Solutions in~\eqref{group.A0} do exist; in
fact, by the uniqueness we conclude that $V_{19}$ is Rohn's configuration.

For the rest of the proof, we subdivide the index set $\IS$ into three
subsets
\[*
\IS_0:=\{1,2,3,4\},\quad\IS_1:=\{5,6,7,8\},\quad\IS_2:=\{9,10,11,12\}
\]
and, for $n=0,1,2$, let
$\one_n:=12\cdot\vv{\IS_n}\in\R^{12}=\bH_{12}\otimes\R$.
This splitting is preserved by all three automorphism groups.
For $V=V_{17}$, we have
\roster*
\item
$\Q V^\perp\cap\SSf\dual=\bA_1(4)$, \ie, $E_0$ is still irreducible,
see~\eqref{eq.Cset};
\item
$\ker M$ is generated by
$\bu:=\one_1-\one_2$.
\endroster
Any solution in the torsion free group~\iref{group.A0**} would have
$p_i=0$, $i\in\IS_0$, violating \eqref{eq.T.distinct}.
All invariant factors of~$M$ divide $6=2\cdot3$ and
the kernel $\ker(M\otimes\FF_3)$
is generated by $\bu$ and
$\one_0-\one_1$, so that
the $\R/\Z$-components of the four points $p_i$, $i\in\IS_0$,
take but two distinct values. Hence, only group~\eqref{group.A0} may
(and does) have solutions satisfying~\eqref{eq.T.distinct}: \eg,
start with Rohn's solution
\[*
\bold{p}\in E_0^{12}=(\R/\Z)^{12}\times(\R/\Z)^{12}
\]
for $V_{19}$ and shift it along the $2$-parametric subgroup
$(\R\bu/\Z)\times(\R\bu/\Z)$.

Finally, for $V=V_{16}$, we have
\roster*
\item
$\Q V^\perp\cap\SSf\dual=\bA_2\dual(4)$, with
$\pm(4\cdot\vv{\IS_n}-\vv{\IS})$ as square $(-8/3)$ vectors;
\item
$\ker M$ is generated by
$\bu:=\one_1-\one_2$ and $\bv:=\one_0-\one_1$.
\endroster
From the first assertion, by~\eqref{eq.Cset},
the cubic~$E_0$ either is irreducible or
splits into a $(-3)$- and $(-4)$-component (a line~$C_0$ and a conic) or three
$(-3)$-components (lines) $C_0$, $C_1$, $C_2$,
so that $p_i\in C_n$ for $i\in\IS_n$ (whenever $C_n$ is present).

As above, we conclude that the projections of $p_i$ to an $\R$-summand or a
finite group $\Z/2$ or $\Z/3$ (the latter selects a component of~$E_0$) are
constant within each set $\IS_n$, ruling out groups~\iref{group.A0**},
\iref{group.A1*}, \iref{group.A2*}.
All invariant factors of~$M$ divide~$2$; hence, for each $n=0,1,2$, the
projections of $p_i$, $i\in\IS_n$, to an $\R/\Z$-factor take but two values.
This rules out groups~\iref{group.A0*},
\iref{group.A1}, \iref{group.A2}. In the
remaining group~\iref{group.A0} a solution is constructed as in the case
$V=V_{17}$.
\endproof

\remark\label{rem.T.brute.force}
The
\ser3-series is the only one where we failed to compute all subsets
$V\subset\VecZ(\SSf,\eta)$ satisfying~\eqref{eq.V}: the counts are huge.
Instead, we have computed, separately, the sets $V$ satisfying~\eqref{eq.V}
and such that
\roster
\item\label{i.>=17}
the cardinality $\ls|V|\ge17$, or
\item\label{i.K4}
the corresponding graph contains $K(4)$, or
\item\label{i.no.K3}
the corresponding graph does not contain~$K(3)$.
\endroster
Recall that
the bound found in \cite[Lemma~4.2]{degt:sextics} is based on a simple observation
that \eqref{eq.V} implies that each
point $\pp\in\IS$ is contained in at most five sets $s\in V$. Hence,
in~\iref{i.>=17} there must be at least three points~$\pp$ contained in exactly
five sets each, and we start with
$14$-element sets~$V$ satisfying this property. We find over $90,000$ sets~$V$
satisfying~\eqref{eq.V} and such that $\ls|V|\ge17$;
five have $\ls|V|=20$.
This indicates that \eqref{eq.V} alone does not reflect the geometry of
quartics very well.

For~\iref{i.K4}, we start with a $4$-element set corresponding to the
graph~$K(4)$. It is due to Lemmata~\ref{lem.K3} and~\ref{lem.T.smooth}
that \iref{i.K4} and \iref{i.no.K3} do complement each other
\emph{as long as we want $\ls|\VZ(\XXf)|=16$ and $\ls|\Fn\XXf|\ge28$,
hence $\Eset(\XXf)=\varnothing$.}
\endremark

%%%%%%%%%%%%%%%%%%%%%%%%%%%%%%%%%%%%%%%%%%%%%%%%%%
%%%%%%%%%%%%%  3.4
%%%%%%%%%%%%%%%%%%%%%%%%%%%%%%%%%%%%%%%%%%%%%%%%%%

\subsection{The $\bX$-series}\label{s.X}

To complete the proof of \autoref{lem.Sigma}, we use $T=\bA_1\oplus\bD_1$, arriving at
the three candidates,
\[*
\SSf=\bE_7\oplus\bA_3,\quad
\bE_8\oplus\bA_1\oplus\bD_1,\quad\text{or}\quad
\bD_9\oplus\bA_1,
\]
with $\bigl|\Vec(\SSf,\Gl)\bigr|=4$, $1$, or $0$, respectively.
By \autoref{lem.K-lines} and \eqref{eq.sets},
we have $\SSf=\bE_7\oplus\bA_3$.

We take for $\lambda$ one of the two generators of $\discr\bA_3=\Z/4$;
then necessarily $\eta=-\lambda+\Ga$, where $\Ga\in\discr\bE_7$ is the generator.
Analyzing the shortest representatives of the discriminant classes, we find
that
\[*
\VecK(\SSf,\lambda)=\VecK(\bA_3,\lambda),\qquad
\VecZ(\SSf,\eta)=\VecK(\bE_7,\Ga)\times\VecK(\bA_3,-\lambda).
\]
Furthermore, viewing $\bE_7$ as $\be^\perp\subset\bE_8$ for a fixed root
$\be\in\bE_8$, we have a bijection
\[*
\VecK(\bE_7,\Ga)\overset\cong\longto
\bigl\{e\in\bE_8\bigm|e^2=-2,\ e\cdot\be=1\bigr\},\quad
l\mapsto l-\be/2.
\]
Together with~$\be\in\bE_8$, the image of a subset $V\subset\VecK(\bE_7,\Ga)$
satisfying~\eqref{eq.V} constitutes a Dynkin diagram~$D$, elliptic or affine,
other than $\tA_1$ and such that
all vertices are adjacent to $\be$. Clearly,
\[
\text{either $D=\tA_2$ or $D\subset\tD_4$};\quad\text{hence},\quad
\bndK(\bE_7,\Ga)=4.
\label{eq.E7}
\]
To the $\bA_3$ summand we apply the machinery of~\autoref{s.T}, with
the new index set $\IS:=\{1,\ldots,4\}$. We have
\[*
\VecK(\bA_{3},\lambda)=
 \bigl\{4\cdot\vv{\{p\}}-\vv\IS\bigm|p\in\IS\bigr\},
\quad
\VecK(\bA_{3},-\lambda)=
 \bigl\{\vv\IS-4\cdot\vv{\{q\}}\bigm|q\in\IS\bigr\},
\]
henceforth regarding both as subsets of~$\IS$.
Then, by \cite[Lemmas~2.6]{degtyarev-90} (\cf. also \autoref{rem.K-lines}),
\roster*
\item
%the set
$\VK(\XXf)$ consists of $r\le4$ points
$p_1<p_2<\ldots<p_r=4$; let $p_0:=0$;
\item
%the set
$\Eset(\XXf)\cap\bA_3$ consists the $(4-r)$ divisors $e_i-e_{i+1}$,
$i\in\IS\sminus\VK(\XXf)$;
\item
%if a $(-2)$-line $s\in\VZ(\XXf)$ contains a point $p\in(p_{k-1},p_k]$,
%it contains the interval $(p_{k-1},p]$.
the projection of $\VZ(\XXf)$ to $\VecK(\bA_3,-\lambda)$ is contained in
the $r$-element set $\{p_0+1,\ldots,p_{r-1}+1\}\subset\VecK(\bA_3,-\lambda)$.
\endroster
Geometrically, the $(-2)$-lines that project to
a point $p_{k-1}+1$ intersect
\emph{in~$\XXf$} (but possibly not in~$\tX$) the $(-1)$-line $p_k$.
The set $\Eset(\XXf)\cap\bA_3$ is a disjoint union of $\bA$-type Dynkin diagrams;
they are what is called the
\emph{essential singularities} in~\cite{degtyarev-90}, \ie, the singular
points of $\XXf(O)$ contained in the tangent plane $C_O\XXf$.

%Together with~\eqref{eq.E7}, this description gives us the following
%immediate consequence.
Invoking~\eqref{eq.E7}, we have the following
immediate consequences.

\lemma\label{lem.X.-1}
If a quartic $\XXf\in\ser4\rat$ is relatively smooth, then $\ls|\FnK\XXf|=4$. More
precisely, $\ls|\FnK\XXf|=4$ if and only if
%$\tX$
$\XXf(O)$
has no singularities on \rom(the preimage of\rom)
the tangent plane $C_O\XXf$.
\done
\endlemma

\lemma[see \cite{gonzalezalonso}]\label{lem.X}
For $\XXf\in\ser4\rat$, one has $\ls|\Fn\XXf|\le5\ls|\FnK\XXf|\le20$.
This bound is sharp.
\done
\endlemma

%\lemma\label{lem.X.tangent}
%Let $\XXf\in\ser4\rat$ and let $\LLL_1,\LLL_2$ be a pair of lines such that
%$\LLL_1\cdot\LLL_2=1$. If $\LLL_1,\LLL_2$ intersect \emph{in $\XXf$} the same
%$(-1)$-line $\LLL$, then the plane spanned by $\LLL_1,\LLL_2$ is tangent to~$\XXf$
%along~$\LLL$.
%\done
%\endlemma

In particular, we conclude that a quartic $\XXf\in\ser4\rat$ with at least
$16$ lines has no exceptional singularities in the sense
of~\cite{degt:quartics2}.
%The following theorem is proved by brute force, see
%\autoref{rem.brute.force}.

\theorem\label{th.X}
If $\XXf\in\ser4\rat$, then one has either $\ls|\Fn\XXf|\le18$
or $\ls|\Fn\XXf|=20$.
%\mnote{\3corrected, updated; doesn't look nice anymore
%\endgraf
%Judging by \autoref{th.T}, must be far from sharp, but I have no idea how
%this can be improved}
Furthermore\rom:
\roster
\item\label{X.12}
if $\ls|\Fn\XXf|\ge16$, then $\ls|\FnK\XXf|=4$\rom;
\item\label{X.14}
if $\ls|\Fn\XXf|=18$, then $\sing(\XXf)$ is as in
\autoref{add-non-K3-projective-double}\iref{add.18}\rom;
%there are at most four $14$-line configurations $\FnZ\XXf$
%\rom(at most six $\Aut_\eta\SSf$-orbits of sets $V$ as in
%\autoref{rem.brute.force}\rom)\rom;
%one has $\Q V=\SSf\otimes\Q$ and $\ls|\Emax|\le1$\rom;
\item\label{X.16}
if $\ls|\Fn\XXf|=20$, then $\sing(\XXf)$ is as in
\autoref{add-non-K3-projective-double}\iref{add.20}.
%there is a unique $16$-line configuration $\FnZ\XXf$
%\rom(at most two $\Aut_\eta\SSf$-orbits of sets $V$ as in
%\autoref{rem.brute.force}\rom)\rom;
%one has $\Q V=e^\perp\subset\SSf\otimes\Q$ and $\Emax=\{\pm e\}$ for a
%certain root $e\in\bE_7$.
\endroster
%In the latter case, $\XXf$ is relatively smooth and
%$\Fn\XXf$ is as follows\rom:
In addition \rom(\cf. \autoref{rem.X} below\rom),
in case~\iref{X.16}, \ie, if $\ls|\Fn\XXf|=20$,
\roster*
%\item
%$\ls|\FnK\XXf|=4$,
\item
each $(-1)$-line intersects four
pairwise disjoint $(-2)$-lines, and
\item
the graph $\FnZ\XXf$ is a generalized quadrangle
$\GQ(3,1)$.
%Furthermore, in this case the set $\VZ$ generates~$\SSf$ over~$\Q$.
%\done
\endroster
\endtheorem

\proof
Statement~\iref{X.12} is given by \autoref{lem.X}. By brute force (see
\autoref{rem.brute.force}), we find two sets
$V\subset\VecZ(\SSf,\eta)$ (a single abstract graph) of size $\ls|V|=16$
and six sets~$V$ (four abstract graphs) of size $\ls|V|=14$. Others have
$\ls|V|\le13$. This immediately implies the assertion that
$\ls|\Fn\XXf|\le18$ or $\ls|\Fn\XXf|=20$.

If $\ls|V|=14$, we have $\Q V=\SSf\otimes\Q$ and $\Emax=\varnothing$ or
$\{e\}$ (for one of the six sets) for a certain root $e\in\bE_7$. Hence,
$\Cset(\XXf)=\varnothing$ and $E_0$ is irreducible, implying that $O$ is of
type $\bX_{1,0}$, $\bX_{1,1}$, or $\bZ^1_{11}$ (see \autoref{tab.resolution}).
In one of the six cases, $\sing(\XXf)$ may also have an extra node.

If $\ls|V|=16$, then $V\cong\GQ(3,1)$ and the description of $\Fn\XXf$ is
given by~\eqref{eq.E7}.
In this case,  $\Q V^\perp$ is spanned by a
certain root $e\in\bE_7$ and  $\Emax=\{\pm e\}$. It follows that either
$\sing(\XXf)$ is as in the previous case or $\Cset(\XXf)=\{\pm e\}$, so that
$E_0$ splits into a $(-4)$- and a $(-2)$-curve, see~\eqref{eq.Cset}, and
$\sing(\XXf)=\bX_{1,2}$.
%
%If the exceptional divisor $E_0$ over~$O$ is reducible, its components~$R_i$
%project to certain vectors $r_i\in\SSf\dual$ orthogonal to $\VZ$: we have
%$r_i^2=-2$, $-3/2$, or~$-2$ for $R_i^2=-2$, $-3$, or $-4$
%(or $R_i\cdot\KK=0$, $1$, or~$2$),
%respectively (see
%\autoref{rem.resolution}).
%For $\ls|\Fn\XXf|=18$ this and \autoref{th.X}\iref{X.14} imply that $E_0$ is
%irreducible, leaving but $\bX_{1,0}$, $\bX_{1,1}$, or $\bZ_{11}$ for the type
%of~$O$, see \autoref{tab.resolution}. Furthermore,
%for one of the four configurations, the set $\sing(\XXf)$ may also have an
%extra node.
%
%If $\ls|\Fn\XXf|=20$, from \autoref{th.X}\iref{X.16} we conclude that $E_0$
%cannot have $(-3)$-\allowbreak components. Hence, either $E_0$ is
%irreducible, leaving room for an extra node, as in the previous case, or
%$E_0$ splits into a $(-4)$-component and exactly one $(-2)$-components: the
%two
%sum up to $-\KK$ and project to $\pm e$; in the latter case,
%\autoref{tab.resolution} implies that $\sing(\XXf)=\bX_{1,2}$.
\endproof

\remark\label{rem.X}
Up to isomorphism, there are two abstract graphs $\Fn\XXf$ as in
\autoref{th.X}: the adjacency to $(-1)$-lines breaks $\FnZ\XXf=\GQ(3,1)$ into
four pairwise disjoint maximal independent subsets, and there are
two $\Aut(\Fn\XXf)$-orbits of such quadruples (which differ by length).
Both are embeddable to~$\SSf$, but we do not know
if both can be realized by lines on a quartic $\XXf\in\ser4\rat$.

Nor do we know if all sets $\sing(\XXf)$ announced may
appear in quartics with $20$ or $18$ lines (\cf. the proof of \autoref{th.T},
where most singularities allowed by lattice theory are eventually ruled out).

\endremark

%%%%%%%%%%%%%%%%%%%%%%%%%%%%%%%%%%%%%%%%%%%%%%%%%%%%%%%%%%
%%%%%%%%%%%%%%%%%  3.5
%%%%%%%%%%%%%%%%%%%%%%%%%%%%%%%%%%%%%%%%%%%%%%%%%%%%%%%%%

\subsection{The $\bJ\*$-series}\label{s.J*}

In this and next (\autoref{s.J} below) cases, the exceptional divisor $E_0$
over~$O$ either is irreducible or has at least one $(-2)$-component, see
\autoref{rem.resolution}.
Hence, the assertion that $\XXf$ is relatively smooth implies, in particular,
that $O$ is of type $\bJ_{2,0}$, $\bJ_{2,1}$, or $\bE_{12}$, see
\autoref{tab.resolution}.

To
complete the proof of \autoref{lem.Sigma}, we use $T=\bA_3$, arriving at
the two candidates,
\[*
\SSf=\bE_8\oplus\bD_1\quad\text{or}\quad
\bD_9,
\]
with $\bigl|\Vec(\SSf,\Gl)\bigr|=1$ or $0$, respectively.
Hence, by \autoref{lem.K-lines} and \eqref{eq.sets},
we have $\SSf=\bE_8\oplus\bD_1$.

We take for $\eta=\lambda$ one of the generators of $\discr\bD_1=\Z/4$.
Then
%set $
\[*
\VK(\XXf)=\VecK(\SSf,\lambda)=\VecK(\bD_1,\lambda)=\{a/4\},
\]
where $a\in\bD_1$ is the generator, and, from analysing the shortest
representatives,
\[*
\VecZ(\SSf,\eta)=
 \bigl(\VecZ(\bE_8,0)\times\VecK(\bD_1,\lambda)\bigr)\cup\VecZ(\bD_1,\lambda).
\]
The last term is $\{-3a/4\}$, and its intersection with $\VZ(\XXf)$, if nonempty,
is (the image of) the special line $\LLX$
intersecting the $(-1)$-line $\LLL_0\in\FnK\XXf$, see \autoref{lem.one.line}.
The subsets $V\subset\VecZ(\bE_8,0)$ satisfying~\eqref{eq.V}
are merely
\emph{parabolic simple graphs}, \ie,
disjoint unions of Dynkin
diagrams, elliptic or affine, other than $\tA_1$. Hence,
\[
\bndZ(\bE_8,0)=12,\quad\text{realized by $4\tA_2\subset\bE_8$}.
\label{eq.E8}
\]
Since $\VecK(\bD_1,\lambda)=\{a/4\}$ is a singleton, we identify the
intersection
\[*
\VZ(\XXf)\cap\bigl(\VecZ(\bE_8,0)\times\VecK(\bD_1,\lambda)\bigr)
\]
with its
projection $V\subset\VecZ(\bE_8,0)$; it is the dual adjacency
graph of the corresponding lines.
Then, due to \autoref{lem.ex}, the set
$\Gamma:=V\cup\Eset(\XXf)$ is also a parabolic simple graph, and we arrive at the
following description of this set:
\roster*
\item
pick a root sublattice
$R=\bigoplus R_i\subset\bE_8$,
%(not necessarily primitive),
where $R_i$ are
the indecomposable components;
\item
let $D_i$ be the Dynkin diagram of $R_i$; convert some of $D_i$ to their
affine counterparts $\tilde{D}_i$;
%\item
%some of the diagrams $D_i$ are to be converted to their affine counterparts
%$\tilde{D}_i\subset V$;
\item
break each component~$D_i$ or $\tilde{D}_i$ into a complementary
pair $D_i'\cup D_i''$ of
induced subgraphs;
\item
let $V=\bigcup_iD'_i$ and $\Eset(\XXf)=\bigcup_iD''_i$, so that $\Gamma$ is
the union of the chosen components $D_i$ or $\tilde{D}_i$.
\endroster
\autoref{lem.K3} imposes the following restrictions:
\roster*
\item
if a line $\LLX$ as in
\autoref{lem.one.line} is present, then
each subgraph $\bA_2\subset V$ is either contained in $\tA_2\subset V$ or
adjacent to a vertex $v\in\Eset(\XXf)$;
\item
if $\LLX$ as in
\autoref{lem.one.line} is not present, then there are no subgraphs
$\tA_2\subset V$.
\endroster
Now, the next lemma is a simple combinatorial exercise.

\lemma\label{th.J*}
For a quartic $\XXf\in\ser5\rat$, one has\rom:
\roster
\item\label{J*.l-cross}
if $\LLX\subset\XXf$ as in \autoref{lem.one.line} is present, then
$\ls|\Fn\XXf|\le12$
or $\ls|\Fn\XXf|=14$\rom;
\item\label{J*.no-l-cross}
if $\LLX\subset\XXf$ as in \autoref{lem.one.line} is not present, then
$\ls|\Fn\XXf|\le11$.
\done
\endroster
\endlemma

\remark\label{rem.J*}
Lattice theoretic techniques do not answer the question
whether the bounds given by \autoref{th.J*} are sharp. In
item~\iref{J*.l-cross},
%if exists, a quartic $\XXf$ with $14$ lines must be
%relatively smooth and the graph $\FnZ\XXf\sminus\LLX=4\tA_2$ is as
%given by~\eqref{eq.E8}.
%%: in
%%addition to $\LLL_0$ and $\LLX$ as in \autoref{lem.one.line}, there are four
%%completely split plane sections in the pencil of planes through~$\LLX$, each
%%containing three extra lines.
%If
%$\ls|\Fn\XXf|=12$, the only
%candidate is $\FnZ\XXf\sminus\LLX=3\tA_2\oplus\bA_1$;
%%possible configuration is obtained from the above
%%by removing two lines from one of the split plane sections.
%this time, $\XXf(O)$ may be smooth or have a node or a cusp.
if $\ls|\Fn\XXf|\ge12$, there are but two candidates for the configuration:
\[*
\FnZ\XXf\sminus\LLX=4\tA_2\qquad\text{or}\qquad
3\tA_2\oplus\bA_1;
\]
in the former case $\XXf(O)$ is smooth, in the latter, it may have a node or
a cusp.
In item~\iref{J*.no-l-cross}, we have but four candidates for the maximal
graph $\FnZ\XXf$:
\[*
\tD_5\oplus\tA_3,\qquad
2\tD_4,\qquad
2\tA_4,\qquad
2\tA_3\oplus2\bA_1.
\]
In the first three cases, $\XXf(O)$ must be smooth; in the last one, it may
have up to two nodes.
We address the realizability question in \autoref{sec-sharp-jstar} below,
see \autoref{th.J*sharp}.
\endremark

%%%%%%%%%%%%%%%%%%%%%%%%%%%%%%%%%%%%%%%%%%%%%%%%%%%%%%%%%%%%%%%%%%%%%%
%%%%%%%%%%%%%%%%%%%%%%%%  3.6
%%%%%%%%%%%%%%%%%%%%%%%%%%%%%%%%%%%%%%%%%%%%%%%%%%%%%%%%%%%%%%%%%%%%%%%

\subsection{The $\bJ$-series}\label{s.J}

To complete the proof of \autoref{lem.Sigma}, we have the same pair of
candidates as in \autoref{s.J*}. According to~\cite{degtyarev-90},
there is a surface $\XXf\in\ser6$ with the set of singularities
$\bJ_{10}\oplus\bD_9$. We conclude that $\SSf=\bD_9$.

We take for $\eta$ one of the generators of $\discr\bD_9=\Z/4$. Since there
are no $(-1)$-lines (see \autoref{lem.one.line}), the other class $\lambda$
makes no sense. The following statement is obtained by brute force (see
\autoref{rem.brute.force}).

\theorem\label{th.J}
If $\XXf\in\ser6$,
then either $\ls|\Fn\XXf|\le13$
or $\ls|\Fn\XXf|=\bndZ(\bD_9,\eta)=16$
and $\Fn\XXf$ is a generalized quadrangle $\GQ(3,1)$.
Furthermore, if $\ls|\Fn\XXf|\ge13$, then $\XXf$ is relatively smooth.
\done
\endtheorem

\remark\label{rem.J}
We do not know
any examples of normal quartics $\XXf\in\ser6$ with many lines: we have one
candidate for $\Fn\XXf$ with $16$ vertices, \viz. $\GQ(3,1)$, and four
candidates with $13$ vertices. Note also that there are three ways to project
$\GQ(3,1)$ to~$\bD_9$: they differ by their stabilizers in $\OG(\bD_9)$.
\endremark

Finding examples of normal quartics $\XXf\in\ser6$ with many lines
is obstructed by
the fact that the
closure of
$\ser6$ contains quartics that are singular along a conic,
as illustrated by the following example.

\example
We assume that the quartic $\XXf$ is given by
%the equations
\eqref{eq-q56} and \eqref{eq-q6} with
$$
Q_2 := (x+y)^2  \quad \mbox{and} \quad H_4 := \frac{1}{4} (x^4 + x^3 y + x y^2 z).
$$
A simple Gr\"{o}bner basis computation shows that $\sing(\XXf) = \mbox{V}(x^2+2 w z)$.
By  \cite[p.~143]{Clebsch-1868} it contains exactly
$16$ lines.
\endexample

\subsection{Irrational quartics}\label{s.irrational}

We
start as in \autoref{s.Sigma} and consider the lattice
$\SSf:=(\Z \hh\oplus\Z\KK)^\perp\subset H_2(\tX)$. For $\XXf$
very general the group
%the lattice
$H_2(\XXf)$ is spanned by $\hh$, the components of $-\KK$,
and the $(-1)$ lines;
hence, $\SSf$ is easily computed.
Indeed, it suffices to observe that,
modulo the radical, the
classes indicated generate a unimodular lattice of the correct rank, see
\autoref{tab.irrational}. The exceptional divisors over $\bX_{2,0}$ and
$\bJ_{4,0}$ are described in \eqref{eq.J4}, and $-\KK=2E+\sum R_i$;
note that $\KK^2=-4$ (resp.\ $-2$) for $\XXf\in\ser4$ (resp.\ $\XXf\in\ser5$).

\remark\label{rem.irrational}
Analyzing the rank of the lattice spanned by
%$\hh$, the components of $-\KK$, and the $(-1)$-lines,
the classes above and referring to \autoref{tab.irrational},
we can make a few geometric conclusions about the
configuration:
\roster*
\item
if $\sing(\XXf)=\bX_{2,0}$, then two $(-1)$-lines intersect~$R_1$ and the two
others intersect $R_2$, see~\eqref{eq.J4};
\item
if $\sing(\XXf)=2\bJ_{10}$, then the two $(-1)$-lines are disjoint.
\endroster
\endremark

The first surprise is that, if $\XXf\in\ser4$, the classes $\hh$ and $\KK$ are
not independent: $\hh=\KK\bmod2H_2(\XXf)$. It follows that $\SSf=\bD_4$, and we
can take for~$\lambda$ any non-zero element of $\discr\bD_4=\Z/2\oplus\Z/2$;
then we have $\ls|\VecK(\bD_4,\lambda)|=8$. For example, interpreting
$\bD_4$ as the maximal even sublattice of $\bH_4$, we can take for $\lambda$
the common discriminant class of the eight square~$1$ vectors $\pm e_i$.

However, in $\SSf$ \emph{there is no room for a class~$\eta$} of
square $-1/4\bmod2\Z$; hence, $\XXf$ has no $(-2)$-lines.
(Alternatively, observe that $\hh\cdot v$ is even for
each $v\in\KK^\perp$.)
Recalling
the relation between $(-1)$-lines and exceptional singularities, \cf.
\cite[Lemma~2.6]{degtyarev-90}, we arrive at the following theorem.

\theorem\label{th.2X}
If $\sing(\XXf)=\bX_{2,0}\oplus\ssing$ or $2\bX_9\oplus\ssing$, then
$\ls|\Fn\XXf|=\ls|\FnK\XXf|$ is as given in \autoref{tab.irrational}.
\done
\endtheorem

If $\XXf\in\ser5$, we have $\SSf=\bA_1\oplus\bD_1$. Let $a$ and $b$ be
generators of the two summands. Then
\[*
\eta=b/4\bmod\SSf,\qquad
\lambda=a/2+b/4\bmod\SSf,
\]
and the lattice has room for two $(-1)$-lines $l_{1,2}:=\pm a/2+b/4$ and
\emph{a single} $(-2)$-line $l\cross:=-3b/4$.
Combining these arguments with \autoref{lem.one.line} and \autoref{ex.J}
below, we arrive at the following statement.

\theorem\label{th.2J}
If $\XXf\in\ser5$ is irrational,
then $\ls|\FnZ\XXf|\le1$ and
%$\ls|\Fn\XXf|$ is bounded
the values of $\ls|\Fn\XXf|$ are
as shown in \autoref{tab.irrational}.
If $\ls|\Fn\XXf|=3$, then the only $(-2)$-line $\LLX$ is the intersection of the
tangent cones at
the two singular points of~$\XXf$.
\done
\endtheorem

%\proof{\3
%The sharpness of the bound for $\sing(\XXf)=2\bJ_{10}$ is given by
%\autoref{ex.J} below, and it remains to prove that the line $\LLX$ does not
%exist if $\sing(\XXf)=\bJ_{4,0}$. A straightforward computation shows that
%$\XXf$ is given by~\eqref{eq-q56}, \eqref{eq-q5} with\mnote{I hope I did this
%right and these are all solutions}
%\[*
%\aligned
%Q&=-3xy,\\
%H&=2x^3z - \frac34x^2y^2
% + h_{040}(4x^2z^2 - 4xy^2z + y^4) - h_{022}(2xz^3 - y^2z^2)
% + h_{004}z^4,
%\endaligned
%\]
%in which case the intersection of~$\XXf$ with the tangent cone
%$C_O\XXf=\{z=0\}$ is
%\[*
%y^2(4h_{040}y^2 + 4wy - 3x^2)=0.
%\]
%The second factor in this equation is always irreducible.
%}
%\endproof

% J_{4,0}
% Q = -3*y*x;
% H = 2*x^3*z - 3/4*x^2*y^2 + 4*x^2*z^2*h[0, 4, 0] - 4*x*y^2*z*h[0, 4, 0] - 2*x*z^3*h[0, 2, 2] + y^4*h[0, 4, 0] + y^2*z^2*h[0, 2, 2] + z^4*h[0, 0, 4];

\example\label{ex.J}
Consider the quartic $\XXf$ given by the equation
\[*
w^2z^2 + wy^3 + x^3z + q_{11}wxyz + h_{220}x^2y^2 = 0.
\]
It is immediate that both the equation itself and the one obtained from it by
the change of variables $w\leftrightarrow z$, $x\leftrightarrow y$ are as
in~\eqref{eq-q56}, \eqref{eq-q5}; hence, $\XXf$ has two singular points of
type~$\bJ_{10}$. The intersection $w=z=0$ of the two tangent cones lies in
$\XXf$ if and only if $h_{220}=0$. Hence, $\XXf$ can have either two or three
lines.
\endexample

% 2J_10
% w^2*z^2*g[0, 0, 2, 2] + w*x*y*z*g[1, 1, 1, 1] + w*y^3*g[0, 3, 0, 1] + x^3*z*g[3, 0, 1, 0] + x^2*y^2*g[2, 2, 0, 0];

\subsection{Proof of \autoref{add-non-K3-projective-double}}\label{proof-add-non-K3-projective-double}
Statement~\iref{add.8}
follows from \autoref{lem.rational} and
\autoref{tab.irrational} (the $\ls|\Fn\XXf|$-column
is given by Theorems~\ref{th.2X} and~\ref{th.2J}).
Statement~\iref{add.16} results from \autoref{tab.bounds} (\eg, the
$E$-column, see \autoref{cor.Elkies}).
All other assertions are given by \autoref{th.X}.
\qed

\newcommand{\sss}{{\mathfrak b}}
\newcommand{\ddd}{{\mathfrak c}}
\newcommand{\qqf}{\ser5\rat}
\def\RRR{R}

%%%%%%%%%%%%%%%%%%%%%%%%%%%%%%%%%%%%%%%%%%
%%%%%%%%%%%%%%%%%   4
%%%%%%%%%%%%%%%%%%%%%%%%%%%%%%%%%%%%%%%%%%

\section{Sharp bound for the \ser5-series \via\ $\sss$-functions}\label{sec-sharp-jstar}

In this section,
we apply the ideas from \cite{Segre,RS-quintics}
to count lines on  $\XXf \in \qqf$.

\theorem \label{th.J*sharp}
A quartic $\XXf \in \ser5$ has at most $12$ lines. This bound is sharp.
\endtheorem

Prior to the proof we collect a few
useful facts.
In view of \autoref{th.J*} and \autoref{th.2J}, we can assume that
$\XXf\in\ser5\rat$ and that it has a line $\LLX\subset C_O\XXf$ as in
\autoref{lem.one.line}, meeting all other lines. Hence, upon rescaling
and by~\eqref{eq.LLX},
\[
h_{301}=1,\qquad h_{220}=0.
\label{eq.h220}
\]
%Recall that all lines on $\XXf \in \qqf$  meet the line $\LLX$ (see \autoref{lem.one.line}).
Moreover, one can easily check that
\begin{equation} \label{eq-smooth along the line}
	\text{the quartic $\XXf$ has no singularities on the line $\LLX$}.
\end{equation}
Consider the morphism
\begin{equation} \label{eq:fibration-deg-3}
	\pi\colon   \XXf \to \PP^1
\end{equation}
given by the linear
system $|{\mathcal O}_{\XXf}(1) - \LLX|$. Its fibers are planar cubics.
We follow \cite{Segre} and say that $\LLX$ is of the \emph{second}
(resp.\ \emph{first}) \emph{kind}  if
%and only if
it is contained in the closure  of the flex
locus of the smooth fibers of
%the fibration
\eqref{eq:fibration-deg-3} (resp.\ otherwise).

The restriction of the fibration \eqref{eq:fibration-deg-3} to the line $\LLX$ defines
the triple cover
\begin{equation} \label{eq-map-deg-d}
	\pi|_{\LLX} \colon \LLX \rightarrow \PP^1.
\end{equation}
By the Hurwitz formula its ramification divisor $\RRR$
has degree
four. One can
%easily
see that the
intersection
point
%$Q_0$  of $\sLLL_0$ and $\LLX$ comes up in $\RRR$ with multiplicity two.
$Q_0:=\sLLL_0\cap\LLX$ has multiplicity~$2$ in~$\RRR$.
Thus, we
have
two possibilities: the support of
%the divisor
$\RRR$  consists of either three points (the ramification type  $(2,1^2)$---see
\cite{RS}) or two points
(the ramification type $(2^2)$):
$$
\text{either \ $\RRR = 2 Q_0 + Q_1 + Q_2$ \ or \ $\RRR = 2 Q_0 + 2Q_1$}.
%\mbox{ where } Q_0 \in \sLLL_0 \cap \LLX \, .
$$
%In particular, one can easily check that if $\XXf$ is given by the equations \eqref{eq-q56} and \eqref{eq-q5} with $h_{301} =1$ (we can make such assumption by rescaling the variables), then the existence of $\LLX$ is equivalent to the condition $h_{220} = 0$. Moreover, if the latter holds,
Assuming~\eqref{eq.h220},
the line $\LLX$ is of ramification type $(2^2)$ if and only if
\begin{equation} \label{eq-ram-type}
	3h_{121} = h_{130} (h_{130} q_{11}^2 - 2 h_{211} q_{11} +3 q_{02}) + 3 h_{040} q_{11} +h_{211}^2 .
\end{equation}

\lemma  \label{lem-firstkind}
If $\LLX$ is a line of the first kind  on $\XXf$, then it is met by
at most $9$  other lines on $\XXf$\rom;
hence, $\ls|\Fn\XXf|\le10$.	
\endlemma

\begin{proof}[Proof {\rm(see \cite[p.~88]{Segre}, \cite[Lemma~5.2]{RS})}]
Clearly, each intersection point $\LLL\cap\LLX$ is in the closure of the flex
locus of the smooth fibers of~\eqref{eq:fibration-deg-3}.
Assuming~\eqref{eq-q56}, \eqref{eq-q5}, and~\eqref{eq.h220}, the resultant of
the restriction to~$\LLX$ of the equation of a fiber
of~\eqref{eq:fibration-deg-3} and its Hessian has degree~$8$ in a
parameter
that equals to~$\infty$ at $C_O\XXf$.
Together with $\LLL_0\subset C_O\XXf$ this makes at most $9$ lines.
%We assume that $\XXf$ is given by the equations \eqref{eq-q56} and \eqref{eq-q5} with $h_{301} =1$ and $h_{220} = 0$.
%	Each line $\ell' \neq \LLX$ that meets $\LLX$ is a component of a fiber of \eqref{eq:fibration-deg-3}. In particular,
%	it meets $\ell$ in a point where both the equation of the cubic curve
%	in question
%	(\ie, the fiber \eqref{eq:fibration-deg-3}) and its hessian vanish.
%	By direct check the resultant of the restrictions of both polynomials to the line $\LLX$ is of degree $8$
%	w.r.t the parameter that parametrizes the planes $\neq  C_O\XXf$  that contain the line $\LLX$ (see \cite[p.~88]{Segre}, \cite[Lemma~5.2]{RS}).
%	Since the plane $C_O\XXf$ contains exactly two lines on $\XXf$ the proof is complete.
\end{proof}

To deal with lines of the second kind we use the rational functions
$\sss_0$, $\sss_1$
%that are
introduced in \cite[Definition~3.3]{RS-quintics} (see also \cite[Remark~3.5]{RS-quintics}).

\lemma  \label{lem-sec211}
If $\LLX$ is a line of the second kind and of
ramification type $(2,1^2)$, then it is met by at most $11$  other lines on $\XXf$. 	
\endlemma

\begin{proof}
Upon the coordinate change $w \mapsto w - h_{130} x - h_{040} y$ in \eqref{eq-q56}
 the ideal of  $\LLX$ is generated by
$x, w$.
Computing
the functions $\sss_0$, $\sss_1$,
we find that their
denominators vanish only at the ramification points of \eqref{eq-map-deg-d}.
Next, we
apply \cite[Proposition~3.9]{RS-quintics} and solve the system of
equations given by the vanishing of the coefficients of the numerator of
$\sss_0$, which has degree~$8$.
%given by the condition
%$\sss_0 = 0$.
%(The numerator of  $\sss_0$ is of degree $8$.)
Substituting the resulting relations between the coefficients of  \eqref{eq-q56} into
$\sss_1$, combined with \cite[Proposition~3.7]{RS-quintics}, shows that at most seven
lines on $\XXf$ meet $\LLX$ away from the support of $\RRR$.
By definition of the ramification type, there are at most three lines on $\XXf$ through
each of the simple
points
$Q_1,Q_2\in\RRR$
(one of them
being $\LLX$ itself) and
exactly two
(\viz. $\LLX$ and  $\sLLL_0$) through the
double
point~$ Q_0=\LLX\cap\sLLL_0$.
This
makes
at most $12$ lines meeting $\LLX$,
and we recall that \emph{exactly}
$12$ lines canot meet $\LLX$ by \autoref{th.J*}.
\end{proof}

\remark \label{rem-211}
An elementary but tedious computation shows that,
under the assumptions of \autoref{lem-sec211}, if three lines on $\XXf$  run
through a reduced point in $\RRR$, say $Q_1$, then at most five lines meet
$\LLX$ away from the support of $\RRR$. In particular, if $\LLX$ is a line of
%the second kind and of
ramification type $(2,1^2)$, then it is met by at
most $10$  other lines on $\XXf$.	
We omit
the details
to keep our
exposition compact.
\endremark

\lemma\label{lem-sec22}
If $\LLX$ is a line of the second kind and of ramification type $(2^2)$,
then it is met by at most $11$  other lines on $\XXf$.
Moreover, if
it is met by exactly $11$ lines,
then exactly one line  $\LLL\neq \LLX$ on $\XXf$ runs through~$Q_1$.
%the point $Q_1$. 	
\endlemma

\begin{proof}
%Again we assume that $\XXf$ is given by the equations \eqref{eq-q56} and \eqref{eq-q5}
%where $h_{301} =1$,  $h_{220} = 0$
%and \eqref{eq-ram-type} holds.
%We change the  variables as in the proof of \autoref{lem-sec211} and
%calculate the functions $\sss_0$, $\sss_1$.
We assume~\eqref{eq-q56}, \eqref{eq-q5}, \eqref{eq.h220},
and~\eqref{eq-ram-type}, change the variables as in the proof of
\autoref{lem-sec211}, and compute $\sss_0$, $\sss_1$.
Solving the system $\sss_0\equiv0$ on $\LLX$,
substituting to~$\sss_1$, and dropping the factors vanishing at $Q_0$
or~$Q_1$ from the numerator,
%We solve the system of equations given by vanishing of $\sss_0$ along the line $\LLX$.
%After substitution of the resulting conditions into $\sss_1$ and
%factoring the forms that vanish in $Q_0$ or $Q_1$  from  the numerator of $\sss_1$
we obtain a degree-$9$ polynomial.
Thus, by
\cite[Proposition~3.7]{RS-quintics}, $\LLX$ is met by at most $10$ lines
(one of them being $\LLL_0$)
away from $Q_1$.

Thus,
it remains to show that
tangent space $T_{Q_1} \XXf$ contains at most two lines on $\XXf$.
%Let us assume that
Otherwise, the quartic curve
$\XXf \cap T_{Q_1} \XXf$ splits into four lines.
By a direct computation
similar to \autoref{rem-211},
the condition that $\sss_0$ vanishes along
%the line
$\LLX$ and the
%cubic residual to $\LLX$ in $\XXf \cap T_{Q_1} \XXf$ consists of
residual cubic in $\XXf \cap T_{Q_1} \XXf$ splits into
three lines implies that $\sss_1$
%vanishes it at most six points
has at most six zeros
away from $Q_0$, $Q_1$.
%Altogether
Thus,
we have at most $(1+6+3)$ lines $\LLL\neq \LLX$ on $\XXf$ that meet $\LLX$.
\end{proof}

\remark
Lemmata~\ref{lem-firstkind}, \ref{lem-sec22}
and \autoref{rem-211} refining \autoref{lem-sec211}
imply that,
%on an $\XXf \in \qqf$ with twelve lines the line
if $\XXf \in \qqf$ has $12$ lines, then
$\LLX$ must be a line of the second type and ramification type $(2^2)$
with exactly two lines on each of the tangent planes $T_{Q_0}$, $T_{Q_1}$.
This resembles the case of smooth quartic surfaces---cf. \cite[Proposition~4.1]{RS}.
%Indeed, such a configuration does appear on $\XXf \in \qqf$, as the following example shows.
\endremark

\example
Consider $\XXf$ given by~\eqref{eq-q56}, \eqref{eq-q5} with $Q_2=(x+y)y$ and
$H_4$ given by
$$
%H_4 :=
-\frac{4}{27} z^4 - \frac{19}{9} y^2 z^2 +  \frac{8}{3} x y^2 z + x^3 z -  \frac{1}{3} y^4 + 2 x^2 z^2 +  \frac{4}{3} x z^3 + 3 x^2 y z -
\frac{4}{3} y z^3  .
$$
One can easily check that \eqref{eq-ram-type} holds
and
the line $\LLX$
given by $z = y-3 w = 0$
is
%a line
of the second kind. In order to see that it is met by exactly
eleven other lines on $\XXf$ one can follow
\latin{verbatim} the approach in \cite[Example~6.3]{RS-quintics}: one checks that
$\sss_1$ has nine
%pairwise different
simple
zeroes away from $Q_0$, $Q_1$
(the  discriminant of its numerator does not vanish), whereas exactly two
lines on $\XXf$ run through $Q_1$ (\cf. the proof of \autoref{lem-sec22}).
\endexample

\proof[Proof of \autoref{th.J*sharp}]
By \autoref{th.J*} we can assume that the
set-theoretic intersection
$C_{O}\XXf\cap\XXf$  consists of two lines.
Lemmata~\ref{lem-firstkind},~\ref{lem-sec211} and \ref{lem-sec22} combined with
\autoref{lem.one.line}\iref{i.all.others}
(and \autoref{th.2J})
complete the proof.
\endproof

%%%%%%%%%%%%%%%%%%%%%%%%%%%%%%%%%%%%%%%%%%%%%%%%%%%%%%%%%%%%%%%%%%%%
%%%%%%%%%%%%%%%%%%%%%%%  5  Quartics with non-isolated singularities
%%%%%%%%%%%%%%%%%%%%%%%%%%%%%%%%%%%%%%%%%%%%%%%%%%%%%%%%%%%%%%%%%%%%

\section{Quartics with non-isolated singularities}\label{S.nonisolated}

In accordance to the general paradigm, and \emph{unlike \autoref{s.ni.line}
below}, by a line we still understand a degree~$1$ curve in $\Cp3$.
There is extensive literature on lines on complex quartic surfaces with
one-dimensional singular locus (see \cite{top-11}, \cite[\S8.6]{dolgachev} and
the bibliography therein).
Below we recall a few basic facts to maintain our exposition self-contained.
To shorten the notation, we adopt the following addendum to
\autoref{conv.types}.

\convention\label{conv.QL}
In addition to \autoref{conv.types}, we say that a quartic $\XXf$
that is  not ruled by lines
is in the
\roster*
\item
\ser2-series, if $\XXf$ has a line~$L$ of double points.
\endroster
\endconvention

%%%%%%%%%%%%%%%%%%%%%%%%%%%%%%%%%%%%%%%%%%%%%%%%%%%%%%%%%%%
%%%%%%%%%%  5.1
%%%%%%%%%%%%%%%%%%%%%%%%%%%%%%%%%%%%%%%%%%%%%%%%%%%%%%%%%%%

\subsection{Taxonomy of non-normal quartics}\label{s.ni.taxonomy}
By Bertini's theorem, a (reduced) curve contained in
the singular locus $\sing(\XXf)$
of a quartic~$\XXf$ has degree at most $3$.
By \cite{top-11},
if $\XXf$ contains (at least) either
\roster*
\item
a twisted cubic of double points, or
\item
a conic and a line of double points, or
\item
a line of triple points, or
\item
two skew lines of singular points,
\endroster
then it
is ruled by lines.
By \cite[p.~176]{top-11}, if
\roster*
\item
$\XXf$ is singular along three concurrent lines,
\endroster
then $\XXf$ is either a cone or
Steiner's Roman surface that contains exactly three lines.
On the other hand, it was shown by Clebsch (\cite[p.~143]{Clebsch-1868}) that if
\roster*
\item
the only one-dimensional component of $\sing(\XXf)$ is a smooth conic
(the so-called \emph{cyclide quartic surface}---see \cite[\S8.6.2]{dolgachev}),
\endroster
then $\XXf$ contains at most $16$ lines (see also  \cite[Lemma~4.3.b]{gonzalezalonso}).
Finally,  if
\roster*
\item
$\XXf$ is singular along two coplanar lines,
\endroster
then it either is ruled by lines (see \cite[\S3.2.6]{top-11}) or
contains at most $18$ lines. The letter claim follows by a
direct determinant computation as in the proof of
\cite[Lemma~3.7]{gonzalezalonso}.
This list exhaust all options but one, and we conclude that
\[
\text{if $\XXf$ is not ruled by lines and
has more than $18$ lines,
%$\ls|\Fn\XXf|>18$,
then $\XXf\in\ser2$}.
\label{eq.ser2}
\]

%\endremark

%To simplify our notation we use the following notation.

\subsection{Lines on
%surfaces
quartics with a line~$L$ of double points}\label{s.ni.line}
Here we study the graph $\Fn\XXf$ for $\XXf\in\ser2$,
with a line~$L$ of double points.
It is important to observe that $L$ itself is not a line in the sense of
\autoref{lem-minustwo}: its pull-back in $\tX$ is an elliptic curve.
In other words,
\begin{equation} \label{eq-L-away-from-Fano}
\text{$L$ does not define a vertex of $\Fn\XXf$}.
\end{equation}	
Thus,
in accordance with
\eqref{eq-fano-on-resolution},
%in all counts in \autoref{s.ni.line}, we do \emph{not} include the line~$L$ itself to
in this section we do \emph{not} include $L$ itself to
$\Fn\XXf$:
for the ultimate statements, including \autoref{thm-27-nonisolated},
an extra $1$
should be added to all counts/bounds.
For the concept of relative smoothness
(\latin{aka} lack of exceptional $(-2)$-divisors),
instead of $\XXf(O)$ we use the
normalization $\XXf(L)$.

Recall that $\ls|\FnK\XXf|\le16$ and a general quartic $\XXf\in\ser2$ has
exactly sixteen $(-1)$-lines (which are those intersecting~$L$),
see \cite[Lemma 3.7]{gonzalezalonso}. These lines appear in pairs $\LLL'$,
$\LLL''$, so that $\LLL'\cdot\LLL''=1$, constituting the singular fibers of
\[
\text{the conic bundle $\XXf\dashrightarrow\Cp1$ given by the projection from~$L$}.
\label{eq.conic.bundle}
\]
Lines from distinct pairs are skew. Generically, each $(-2)$-line intersects
exactly one $(-1)$-line from each pair.

A
%simple
computation using~\eqref{eq.conic.bundle}
%the conic bundle
shows that, if $\XXf(L)$ has simple singularities only,
then
\[*
\hh^2=4,\quad
\KK^2=0,\quad
\hh\cdot\KK=-2,\quad\text{hence}\quad
b_2(\tX)=10.
\]
Consider a general surface~$\XXf$ with sixteen $(-1)$-lines~$\LLL_i$
and at least
one $(-2)$-line $m$. The classes of $h$, $\KK$, $\LLL_i$, and $m$, modulo
radical, generate a unimodular lattice of rank~$10$, which therefore has to
be $H_2(\tX)$. The lattice
\[*
\SSf:=(\Z \hh\oplus\Z\KK)^\perp=\bD_8
\]
is computed directly.
The $(-2)$- and $(-1)$-lines project to vectors of square
\[*
\qZ:=-2\qquad\text{or}\qquad\qK:=-1,
\]
respectively.
Intersections of the projections of distinct
$(-2)$-lines take values $\qZ+2$ or $\qZ+3$, as in \autoref{Sboundslattices}; those of a
$(-2)$-line and a $(-1)$-line take values $\pm1/2$.
Thus, we have an immediate bound based on Elkies~\cite{elkies} (\cf. the
proof of \autoref{cor.Elkies});
it turns out better than \cite[Lemma 3.9]{gonzalezalonso} but worse than
\cite[Example 3.10]{gonzalezalonso}.

%A simple computation shows that, if $\XXf(L)$ has simple singularities only,
%then\mnote{\3to be proved}
%\[*
%\hh^2=4,\quad
%\KK^2=-1,\quad
%\hh\cdot\KK=-2,\quad\text{hence}\quad
%b_2(\tX)=11.
%\]
%Consider a general surface~$\XXf$ with sixteen $(-1)$-lines~$\LLL_i$
%and at least
%one $(-2)$-line $m$. The classes of $h$, $\KK$, $\LLL_i$, and $m$, modulo
%radical, generate a unimodular lattice of rank~$11$, which therefore has to
%be $H_2(\tX)$. The lattice $\SSf$ is computed directly: its maximal root
%system is $\bD_8$ and $\SSf$ is the only, up to isomorphism, index~$2$
%extension
%\[
%\SSf\supset\bD_8\oplus\Z a,\quad a^2=-8,\qquad\text{with}\quad\discr\SSf=\Z/8;
%\label{eq.Sigma.L}
%\]
%the $\bD_8$ summand is canonically defined as the maximal root sublattice
%of~$\SSf$.
%(Alternatively, $\SSf$ can be represented as an index~$2$
%sublattice of $\bE_8\oplus\bA_1$.)
%The $(-2)$- and $(-1)$-lines project to vectors of square, respectively,
%\[*
%\qZ:=-17/8\qquad\text{or}\qquad\qK:=-9/8.
%\]
%Intersections of the projections of distinct
%$(-2)$-lines take values $\qZ+2$ or $\qZ+3$, as above; those of a
%$(-2)$-line and a $(-1)$-line take values $-3/8$ or $5/8$.
%These data let us state an immediate bound based on Elkies~\cite{elkies} (\cf. the
%proof of \autoref{cor.Elkies});
%it turns out better than \cite[Lemma 3.9]{gonzalezalonso} but worse than
%\cite[Example 3.10]{gonzalezalonso}.

\lemma\label{lem.Elkies.L}
For a quartic $\XXf\in\ser2$,
Elkies' bound~{\rm \cite{elkies}}
is $\ls|\FnZ\XXf|\le12$\rom; hence $\ls|\Fn\XXf|\le28$.
\done
\endlemma

%For the lattice theoretic approach,
Next, we proceed as in \autoref{s.Sigma},
taking for $\eta\ne0$ (resp.\ $\lambda$) any square $0$
(resp.\ square $1$)
generator of $\discr\SSf=\Z/2\oplus\Z/2$.
We have
%\[*
%\VecZ(\SSf,0)=\VecZ(\bD_8,0)\qquad\text{and}\qquad
%\bigl|\VecK(\SSf,\lambda)\bigr|=16,
%\]
\[*
\VecK(\bD_8,\lambda)=\bigl\{\pm e_i\bigm|
 \text{$e_i$ is a standard generator of $\bH_8\supset\bD_8$}\bigr\},
\]
accommodating for exactly sixteen $(-1)$-lines. Furthermore, for
each $e\in\VecZ(\SSf,0)$ there are two vectors $v\in\VecK(\SSf,\lambda)$ with
$e\cdot v<0$. In view of \autoref{rem.K-lines}, we have the following
criterion.

\lemma\label{lem.L.smooth}
An
%quartic
$\XXf\in\ser2$ is relatively smooth if and only if
$\ls|\FnK\XXf|=16$.
\done
\endlemma

Next, we compute
\[*
\bigl|\VecZ(\bD_8,\eta)\bigr|=128;
\]
the elements of this set are \emph{some} (not all) roots
in the index~$2$ extension $\bE_8\supset\bD_8$ by $\eta$.
It follows that, as in
\autoref{s.J*}, both $\VZ(\XXf)$ and $\VZ(\XXf)\cup\Eset(\XXf)$
are
parabolic simple graphs, although not any graph may appear.
Using brute force, we find that
\[*
\bndZ(\bD_8,\eta)=10
\]
and there are but three $\OG_\eta(\bD_8)$-orbits of sets
satisfying~\eqref{eq.V}; their graphs are
\[*
2\tD_4,\qquad
\tD_5\oplus\tA_3,\qquad\text{or}\qquad
2\tA_3\oplus2\bA_1,
\]
each implying $\Eset(\XXf)=\varnothing$ by \autoref{lem.ex}.
Thus, taking into account \cite[Example 3.10]{gonzalezalonso},
we have the following ultimate statement, improving \autoref{lem.Elkies.L} by
two more units.

\theorem\label{th.L}
For a quartic $\XXf\in\ser2$, one has
$\ls|\FnZ\XXf|\le10$ and $\ls|\Fn\XXf|\le26$. Both bounds are
sharp. If $\ls|\Fn\XXf|=26$, then the quartic is relatively smooth and
$\FnZ\XXf=2\tD_4$, $\tD_5\oplus\tA_3$, or
$2\tA_3\oplus2\bA_1$.
\done
\endtheorem

\remark\label{rem.L}
Assuming that $\ls|\FnK\XXf|=16$, \ie, $\VK(\XXf)=\VecK(\bD_8,\lambda)$, one
can easily recover the full graph $\Fn\XXf$. One geometric
restriction in the spirit of
\autoref{lem.K3} is that there should be no triangles $K(3)\subset\Fn\XXf$.
We omit details.

We do not know which of the three configurations of lines can be realized: in
the example in~\cite{gonzalezalonso}, only the \emph{number} of lines is
known; it is found by counting the roots of a certain polynomial.
\endremark

%{\2{\bf To do:} to check the anticipated intersections above. To say the usual BS
%about $\XXf(L)$ having simple singularities only or to see what happens if a
%non-simple point is present. To rewrite the whole thing in a nice way.
%}

\subsection{Proof of \autoref{thm-27-nonisolated}}\label{proof-27-nonisolated}
By \eqref{eq.ser2}, we can assume that $\XXf \in \ser2$.
Then, \autoref{th.L}
with an extra $1$ added to the count due to~\eqref{eq-L-away-from-Fano}
implies that $\XXf$ contains at most $27$ lines.
%Indeed, each line on $\XXf$, different from $L$, defines an element of
%$\Fn\XXf$ and we obtain the upper bound of at most $27$ lines on
%$\XXf \in \ser2$ (recall \eqref{eq-L-away-from-Fano}).
Finally,
if the bound is attained, $\ls|\Fn\XXf|=26$, from \autoref{th.L} again we
conclude that $\XXf$
has no singularities away from
the rational curve $L$.
\qed

\remark \label{rem-elkies-rules}
Observe that the non-sharp bound of at most $29$ lines on a quartic
$\XXf \in \ser2$ that we obtain without \GAP~\cite{GAP4.13}
(see \autoref{lem.Elkies.L}), is strong enough for the proof of
the upper bound in
\autoref{thm-non-K3-projective}.
\endremark

\vspace*{1ex}
\noindent
Alex Degtyarev\\
Department of Mathematics, Bilkent University, 06800 Ankara, Turkey\\
degt@fen.bilkent.edu.tr

\vspace*{1ex}
\noindent
S{\l}awomir Rams \\
Faculty of Mathematics and Computer Science, Jagiellonian University,
 ul. {\L}ojasiewicza 6,  30-348 Krak\'{o}w, Poland \\
slawomir.rams@uj.edu.pl

% \vspace*{5ex}
% \noindent
%\filename

\end{document}